\documentclass[3p]{elsarticle}

\usepackage{lineno,hyperref}
\usepackage{siunitx}
\usepackage{booktabs}
\usepackage{colortbl}
\usepackage{tikz}   
\usepackage{float}  
\usepackage{placeins}
\usepackage{graphicx}
\usepackage[caption=false]{subfig}
\usepackage{float}
\usepackage{multirow}
\usepackage{amsmath,mathtools}
\usepackage{amsfonts}
\usepackage[title]{appendix}
\usepackage{soul}
\usepackage{comment}
\usepackage{enumitem}
\usepackage{mathrsfs}

\usepackage{algorithm,algorithmic}

\usepackage{cases}

\hypersetup{pdfauthor={Emmanouil G. kakouris}}

\usepackage[most]{tcolorbox}
\newtcolorbox{highlighted}{colback=yellow,coltext=red,breakable}

\DeclareMathOperator{\Tr}{Tr}
\DeclareMathOperator\supp{supp}

\newtheorem{Remark}{Remark}
\newcommand{\FIG}[1]{Fig.~\ref{#1}}

\newcommand{\EQ}[1]{Eq.~\eqref{#1}}
\newcommand{\EQS}[1]{Eqs.~\eqref{#1}}

\newlength\myindent
\newcommand\bindent{%
  \begingroup
  \setlength{\itemindent}{\myindent}
  \addtolength{\algorithmicindent}{\myindent}
}
\newcommand\eindent{\endgroup}
\newcommand {\R}{\mathbb{R}}

\modulolinenumbers[5]

\begin{document}

\begin{frontmatter}


\title{A high-fidelity material point method for frictional contact problems}


\author[address1]{Emmanouil G. Kakouris \corref{mycorrespondingauthor}}
\cortext[mycorrespondingauthor]{Corresponding author}
\ead{Emmanouil.Kakouris@warwick.ac.uk}

\author[address2]{Manolis N. Chatzis}

\author[address3]{Savvas P. Triantafyllou}

\address[address1]{School of Engineering, University of Warwick, Coventry, CV4 7AL, UK}

\address[address2]{Department of Engineering Science, University of Oxford, Oxford, OX1 3PA, UK}

\address[address3]{Institute for Structural Analysis and Aseismic Research, School of Civil Engineering, National Technical University of Athens, Greece}


\begin{abstract}
A novel Material Point Method (MPM) is introduced for addressing frictional contact problems. In contrast to the standard multi-velocity field approach, this method employs a penalty method to evaluate contact forces at the discretised boundaries of their respective physical domains. This enhances simulation fidelity by accurately considering the deformability of the contact surface, preventing fictitious gaps between bodies in contact. Additionally, the method utilises the Extended B-Splines (EBSs) domain approximation, providing two key advantages. First, EBSs robustly mitigate grid cell-crossing errors by offering continuous gradients of the basis functions on the interface between adjacent grid cells. Second, numerical integration errors are minimised, even with small physical domains in occupied grid cells. The proposed method's robustness and accuracy are evaluated through benchmarks, including comparisons with analytical solutions, other MPM-based contact algorithms, and experimental observations from the literature. Notably, the method demonstrates effective mitigation of stress errors inherent in contact simulations.
\end{abstract}


\begin{keyword}
\texttt{material point method, frictional contact mechanics, extended b-splines}
\end{keyword}

\end{frontmatter}



\section{Introduction} \label{sec:intro}

Robust numerical modelling of frictional contact problems is crucial for various engineering applications, including tyre-road interaction \cite{Yang2020TheValidation}, soil-structure interaction \cite{Zhang2019AnAnalysis,Chatzis2012ModelingProblem}, and metal forming \cite{Marques1990ThreedimensionalForming} amongst many. Despite significant progress in improving the accuracy of numerical simulations in contact mechanics, capturing the kinematic response of solid bodies during contact or impact remains a challenging task \cite{Sauer2013LocalScheme,Duong2019AMethod}.

Errors in stress predictions become more pronounced when dealing with deformable bodies of complex geometries in contact or impact scenarios. Mesh-based simulation methods (e.g., \cite{Laursen2012MortarFormulations,Scolaro2022DevelopmentSolvers,Xing2018AProblems,Aldakheel2020CurvilinearMechanics}) can treat frictional contact, even in large deformations, when the mesh distortion error is not severe \cite{Popp2009AStrategy,DeLorenzis2012AMethod}. However, the fidelity of the simulation rapidly deteriorates due to mesh distortion. To address this, various numerical remedies have been proposed, such as re-meshing, though requiring additional computational time \cite{PavanaChand1998RemeshingProblems}.

The Material Point Method (MPM) \cite{Sulsky1994AMaterials} is a robust approach for solving contact problems by discretising the continuum into Lagrangian material points. These points carry properties such as mass, position, velocity, and stress. Unlike purely meshless methods, MPM employs a background Eulerian grid (computational mesh) for calculating forces and momentum. With its Eulerian-Lagrangian formulation, MPM delivers high-fidelity solutions for domains undergoing large deformations \cite{Wang2019OnDistortion,Coombs2020OnElasto-plasticity}. Notably, MPM excels in contact mechanics problems, enabling the straightforward detection of contact surfaces \cite{Bardenhagen2001AnMaterial, Homel2017FieldgradientMethod, Kakouris2019Phase-FieldEnergy, Moutsanidis2019ModelingField, Nairn2020NewFiltering, Xiao2021DP-MPM:Fragmentation, deVaucorbeil2021ModellingMethod, Guilkey2021AMethod}. However, the MPM has been reported to face three main challenges, which yield severe stress errors (“noise”) at the contact surface, i.e., errors arising from premature contact detection (early contact), grid cell-crossing errors, and errors arising in stresses due to incomplete numerical integration over boundary grid cells when these are occupied by a small physical domain. In this work, we introduce a novel MPM variant that leverages the merits of the Extended B-Spline MPM \cite{Yamaguchi2021ExtendedMethod} and reformulates it within a penalty function contact setting inspired by \cite{Hamad2017AMPM} to simultaneously treat these errors. 

\emph{Early contact errors} refer to the premature contact stresses arising between material points, i.e., before the physical domains actually come into contact, see, e.g., \cite{Huang2011ContactSimulation}. This aspect has been treated using penalty methods, starting from the seminal work of  \cite{Hamad2017AMPM}. In  \cite{Liu2020ILS-MPM:Particles}, the penalty method has been combined with a level-set representation of the boundary surfaces to increase accuracy. Very recently, \cite{Guilkey2021AMethod}, see also \cite{Guilkey2023CohesiveMethod} proposed a penalty contact method where the physical boundary is  discretised into line segments (in 2D) or triangulated surfaces (in 3D). The contact forces are then projected from the discretised contact area back to the Eulerian grid where the governing equations are solved. Although, this implementation provides more accurate estimates of the contact problem, the domain has to be discretised sufficiently fine to approximate a solution region accurately. Very recently, \cite{Chen2023DEM-enrichedMethod}  introduced an MPM method that resides on a Discrete Element Method (DEM) formulation to treat multi-material interactions. 

An alternative approach to penalty methods has been proposed in \cite{Nairn2020NewFiltering} where logistic regression was employed to estimate the contact surface normal vector and to define the domain boundaries. However, the robustness of the logistic regression method remains an open issue. As already explained in \cite{Guilkey2021AMethod}, logistic regression may struggle in scenarios featuring small deformations between stiff objects. In such cases, including scenarios with curved surfaces, the method, along with its predecessors, tends to overestimate the contact area, consequently leading to an underestimation of stress in contact regions. In this work, the solid body boundaries are represented by introducing additional material points, termed boundary material points. Those boundary material points have both mass and volume contribution to the solid body they belong, and are treated in a similar manner to  the materials points in the bulk. Thus, no special treatment is required to update their properties, i.e. position, velocity, etc. Contact forces arise when a material point along the boundary of one solid body penetrates a segment defined by material points along the boundary of another solid body.

\emph{Grid cell-crossing errors} occur when material points cross grid cell boundaries due to non-continuous gradients of the basis functions. To this point, several methods have been introduced to treat this, starting from the Generalised Interpolation Material Point (GIMP), see, e.g., \cite{Bardenhagen2004TheMethod, Charlton2017IGIMP:Deformations}, the Convected Material Point Domain Interpolation (CPDI), see, e.g., \cite{Sadeghirad2011ADeformations,Tran2020AMethod}, the Second-order CPDI (CPDI2), see, e.g, \cite{Sadeghirad2013Second-orderInterfaces}, the Total Lagrangian Material Point Method (TLMPM) \cite{deVaucorbeil2021ModellingMethod}, the $\mathcal{C}_0$ enhancement scheme \cite{Wilson2021DistillationRemedy}, and the Staggered Grid MPM (SGMP) \cite{Liang2019AnMethod}. The CPDI is one of the most accurate variants of the GIMP to date. In CPDI, the shape functions of the background grid are replaced with shape functions defined within each material point's domain, typically assumed to be a parallelogram whose edges require tracking. Due to this issue, the CPDI adds a level of complexity to the MPM as the position of the domain has to be tracked in each computational cycle, is not easily ammenable to a parallel implementation \cite{Homel2017FieldgradientMethod, Nguyen2023TheApplications}  and suffers from mesh distortion errors \cite{Wang2019OnDistortion,Wang2024AnProblems}.

In comparison, B-Spline basis functions are highly effective in mitigating the cell-crossing instability, without the requirement of additional algorithmic manipulations. B-splines are directly employed at the background grid cells and resolve grid cell crossing by providing increased continuity at the interface. The merits of using a higher-order continuity basis in the MPM have been extensively examined in the literature, see, e.g., \cite{Bing2019B-splineMethod, Moutsanidis2020IGA-MPM:Method}. In this work, we adopt this remedy because of its favourable characteristics. In this work, we employ an Extended B-Splines interpolation scheme (EBS) that in addition to providing higher order continuity, can minimise numerical integration errors, even when boundary grid cells are occupied by small physical domains.

\emph{Numerical integration errors} refer to stress errors arising from the incomplete integration over boundary grid cells when occupied by a small physical domain. This occurs when material points occupy a region much smaller than the grid cell. Mitigating this involves techniques such as the ``cut-off" method proposed by \cite{Sulsky1995ApplicationMechanics}, which introduces a model parameter, but it imposes undesirable constraints. Other approaches include the momentum formulation of the material point method \cite{Sulsky1995ApplicationMechanics} and the Modified Update Stress Last (MUSL) approach \cite{Nguyen2023TheApplications}, both adopted in this work (detailed in Section \ref{sec:mpm_implementation}). 

The EBS has been first introduced in \cite{Hollig2002WeightedProblems} within an MPM setting and has been shown to alleviate integration errors, although not within a contact mechanics setting. This work, addresses this issue by originally introducing the EBS in a penalty contact numerical framework. Furthermore, our approach introduces boundary material points to precisely resolve the deformable domain's physical boundary and the evolving contact stresses. Contrary to \cite{Guilkey2021AMethod}, in our work the boundary material points are physical points, have mass. They hence contribute to the stress distribution within the domain. Our benchmark tests verify the algorithm's effectiveness, comparing numerical results with analytical solutions and state-of-the-art MPM contact implementations from the literature.

The manuscript is organised as follows: Section \ref{sec:gov_eqns_frict_cont_probs} introduces the governing equations for frictional contact problems. Subsequently, Section \ref{sec:mpm_implementation} delves into the EBSs-based MPM discretisation scheme and its numerical implementation. Section \ref{sec:numerical_examples} presents numerical examples to verify and validate the proposed method, comparing its efficiency with analytical solutions, other contact algorithms and experimental observations. Finally, Section \ref{sec:conclusions} offers concluding remarks and outlines future directions.

\section{Governing equations for frictional contact} \label{sec:gov_eqns_frict_cont_probs}

\subsection{Principal balance of energy} \label{sec:princ_balance_eng}

\begin{figure}[htbp]
	\centering
	\begin{tabular}{cc}
		\subfloat[\label{fig:gov_eqns_frict_cont_probs:continuum_bodies_mpm_approximation:a}]{
			\includegraphics[scale=1.0]{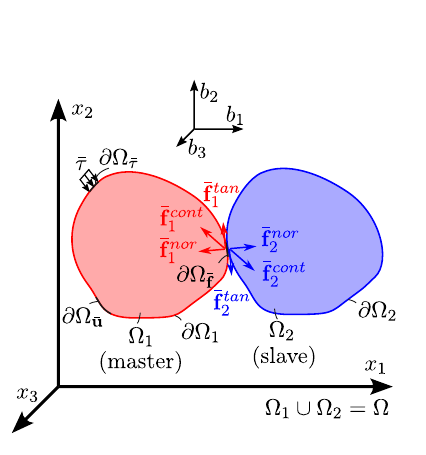}} &
		\subfloat[\label{fig:gov_eqns_frict_cont_probs:continuum_bodies_mpm_approximation:b}]{
			\includegraphics[scale=1.0]{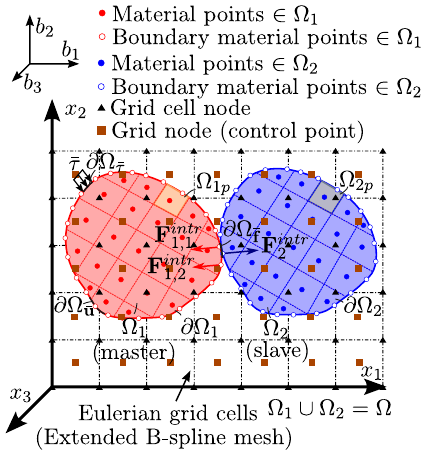}}
	\end{tabular}
	\caption[]{\subref{fig:gov_eqns_frict_cont_probs:continuum_bodies_mpm_approximation:a} Continuum bodies (discrete fields) into contact and \subref{fig:gov_eqns_frict_cont_probs:continuum_bodies_mpm_approximation:b} Material point method approximation.}
    \label{fig:gov_eqns_frict_cont_probs:continuum_bodies_mpm_approximation}
\end{figure}

This section briefly discusses the case of two bodies in contact, as illustrated in \FIG{fig:gov_eqns_frict_cont_probs:continuum_bodies_mpm_approximation}. The method is naturally extendable to handle multiple domains in contact, as demonstrated in Section \ref{sec:numerical_examples}. The study focuses on a deformable domain $\Omega$ with an outer boundary $\partial \Omega = \partial \Omega_1 \cup \partial \Omega_2$, comprising two bodies denoted as discrete fields, namely $\Omega_1$ (master) and $\Omega_2$ (slave). Here, $\Omega = \Omega_1 \cup \Omega_2 \subset \R^{d}$ with $d \in {1, 2, 3}$. At time $t$, the deformable bodies are assumed to be in contact along a contact surface $\partial \Omega_{\bar{\mathbf{f}}}$, as shown in \FIG{fig:gov_eqns_frict_cont_probs:continuum_bodies_mpm_approximation:a}. The normal and tangential contact vectors are denoted by $\bar{\mathbf{f}}_1^{nor} = - \bar{\mathbf{f}}_2^{nor}$ and $\bar{\mathbf{f}}_1^{tan} = - \bar{\mathbf{f}}_2^{tan}$, respectively, with their magnitudes represented as $\bar{\mathbf{f}}_1^{cont} = - \bar{\mathbf{f}}_2^{cont}$.

Assuming contact violation, where one discrete field penetrates another in the region $\partial \Omega_{\bar{\mathbf{f}}}$, the total energy of domain $\Omega$ is penalised using a proportional penalty governed by the function $\mathcal{P}$. This penalty function is expressed as:

\begin{equation}
\label{eqn:penalty_func}
\begin{aligned}
\mathcal{P} \left( \mathbf{u} \right) = \frac{1}{2} \omega^{nor} \int_{\partial\Omega_{\mathbf{\bar{f}}}}\left( \textsl{g}^{nor} \right)^2 d \partial\Omega_{\mathbf{\bar{f}}} + \frac{1}{2} \omega^{tan} \int_{\partial\Omega_{\mathbf{\bar{f}}}}\left(\textsl{g}^{tan}\right)^2 d \partial\Omega_{\mathbf{\bar{f}}}
\end{aligned},
\end{equation}

Here, $\omega$ is a large penalty parameter, $\textsl{g}^{nor}$ is the normal gap function, and $\textsl{g}^{tan}$ is the tangential slip function. The superscripts $\left( \cdot \right)^{nor} $ and $\left( \cdot\right)^{tan}$ denote the normal and tangential directions, respectively.

The normal gap and tangential slip functions, following \cite{Hamad2017AMPM}, depend on the relative distance between a slave point, $s$, on $\partial \Omega_{2\bar{\mathbf{f}}} \subset \partial \Omega_2$, and its projection onto the master segment of $\partial \Omega_{1\bar{\mathbf{f}}} \subset \partial \Omega_1$. The symbols $\partial \Omega_{1\bar{\mathbf{f}}}$ and $\partial \Omega_{2\bar{\mathbf{f}}}$ represent the contact surfaces of $\Omega_1$ and $\Omega_2$, respectively.

\begin{Remark} The formulation in \EQ{eqn:penalty_func} quantifies the extent of constraint violation by employing a penalty coefficient that scales with the constraint error. Detailed information about the gap and slip functions used in this context is available in Section \ref{sec:definition_gap_function}.
\end{Remark}

It follows from \EQ{eqn:penalty_func} that the rate of the penalty function can be written as
\begin{equation}
\label{eqn:rate_penalty_func}
\begin{aligned}
\dot{\mathcal{P}} \left( \dot{\mathbf{u}} \right) = \omega^{nor} \int_{\partial\Omega_{\mathbf{\bar{f}}}}\left( \dot{\textsl{g}}^{nor} \right) d \partial\Omega_{\mathbf{\bar{f}}} + \omega^{tan} \int_{\partial\Omega_{\mathbf{\bar{f}}}}\left(\dot{\textsl{g}}^{tan}\right)d \partial\Omega_{\mathbf{\bar{f}}}.
\end{aligned}
\end{equation}

Hence, the total energy of the system is expressed as 
\begin{equation}
\label{eqn:potential_eng}
\begin{aligned}
\mathcal{T} \left( \mathbf{u} \right) = \mathcal{K} \left( \mathbf{u} \right) + \mathcal{W}^{int} \left( \mathbf{u} \right) - \mathcal{W}^{ext} \left( \mathbf{u} \right) + \mathcal{P} \left( \mathbf{u} \right),
\end{aligned}
\end{equation}
where $\mathcal{K}$, $\mathcal{W}^{int}$ and $\mathcal{W}^{ext}$ are the total kinetic energy, the internal work, and the external to the system work, respectively \cite{Wriggers2006ComputationalMechanics}. The symbol $\mathbf{u} = \mathbf{u} \left( \mathbf{x}, t \right)$ denotes the displacement of a point $\mathbf{x} = \{ x_1, x_2, x_3 \}$ at time $t$. Employing the principle of balance of energy, the strong form for frictional contact is derived as
\begin{equation}
\label{eqn:princ_balance_eng}
\begin{aligned}
\dot{\mathcal{T}} \left( \dot{\mathbf{u}} \right) = \dot{\mathcal{K}} \left( \dot{\mathbf{u}} \right) + \dot{\mathcal{W}}^{int} \left( \dot{\mathbf{u}} \right) - \dot{\mathcal{W}}^{ext} \left( \dot{\mathbf{u}} \right) + \dot{\mathcal{P}} \left( \dot{\mathbf{u}} \right) = 0,
\end{aligned}
\end{equation}
where the symbol $\dot{\left( \cdot \right)}$ denotes differentiation with respect to time.

\subsection{Strong form} \label{sec:strong_form}

Applying the divergence theorem in \EQ{eqn:princ_balance_eng}, performing the necessary algebraic manipulations, and finally
considering that the resulting expression must hold for arbitrary values of $\dot{\mathbf{u}}$, the strong form for frictional contact problems is derived as
\begin{equation}
\label{eqn:strong_form}
\begin{aligned}
\nabla \cdot \boldsymbol {\sigma}  + \mathbf{b} = \rho \ddot{\mathbf{u}} \quad \quad \text{on } \Omega
\end{aligned}
\end{equation}
where $\boldsymbol{\sigma} = {\partial \psi_{el}} / {\partial \boldsymbol{\varepsilon}}$ is the stress field and $\psi_{el}(\boldsymbol{\varepsilon})$ is the elastic energy density. In this work, linear elastic materials are considered and, hence the strain tensor, $\boldsymbol{\varepsilon}$ is expressed as 
\begin{equation}
\begin{aligned}
\boldsymbol \varepsilon = \frac{1}{2} \Big(\nabla \mathbf{u} + \nabla \mathbf{u}^{T} \Big)
\label{eqn:SmallStrainsApprox}
\end{aligned}
\end{equation}
according to the infinitesimal strain theory \cite{Bonet2008NonlinearAnalysis}; the $\left(\nabla\right)$ stands for the gradient operator. In \EQ{eqn:strong_form}, $\mathbf{b} = \{  b_1, b_2, b_3 \}$ are the body forces, $\rho$ the mass density and $\ddot{\mathbf{u}} = d \mathbf{u} / dt$ the acceleration field.

\EQ{eqn:strong_form} is also subjected to the set of boundary and initial conditions defined in \EQ{eqn:initial_conditions}
\begin{equation}
\label{eqn:initial_conditions}
\begin{aligned}
 \begin{cases}
\boldsymbol{\sigma} \cdot \mathbf{n} = \bar{\mathbf{\tau}}, & \text{on } \partial\Omega_{\bar{\mathbf{\tau}}} \\
\mathbf{u} = \bar{\mathbf{u}}, & \text{on } \partial\Omega_{\bar{\mathbf{u}}} \\
\mathbf{u}=\prescript{(0)}{}{  \mathbf{u}   }^{}_{}, & \text{on } \prescript{(0)}{}{  \Omega   }^{}_{} \\
\dot{\mathbf{u}}=\prescript{(0)}{}{  \dot{\mathbf{u}}   }^{}_{}, & \text{on } \prescript{(0)}{}{  \Omega   }^{}_{} \\
\ddot{\mathbf{u}}=\prescript{(0)}{}{  \ddot{\mathbf{u}}   }^{}_{}, & \text{on } \prescript{(0)}{}{  \Omega   }^{}_{}
\end{cases}
\end{aligned}
\end{equation}
where $\mathbf{n}$ is the outward unit normal vector of the boundary, $\bar{\mathbf{u}}$ are the prescribed displacements on $\partial \Omega_{\bar{\mathbf{u}}}$, and $\bar{\mathbf{\tau}}$ corresponds to the set of traction forces at $\partial \Omega \bar{\mathbf{\tau}}$. The notation $\prescript{(t)}{}{ \left( \cdot \right) }^{}_{}$ indicates the quantity $\left( \cdot \right)$ at time $t$.

Furthermore, the strong form, \EQ{eqn:strong_form} is subjected to the kinematic constraints presented in \EQ{CollinearityNormalCond} to \EQ{ComplementaryNormalCond} and \EQ{CollinearityTangCond} to \EQ{ComplementaryTangCond} at the contact surface $\partial\Omega_{\bar{\mathbf{f}}}$ \cite{Wriggers2006ComputationalMechanics}. 

The kinematic constraints of \EQ{CollinearityNormalCond} to \EQ{ComplementaryNormalCond} correspond to the normal contact laws
\begin{subnumcases}{}
\mathbf{e}^{nor}_{1} = - \mathbf{e}^{nor}_{2}, & Collinearity, \space \space \space \space \space \space \space \space on $\partial\Omega_{\bar{\mathbf{f}}}$  \label{CollinearityNormalCond} \\
\bar{\mathbf{f}}_{1}^{nor} = - \bar{\mathbf{f}}_{2}^{nor}, & Collinearity, \space \space \space \space \space \space \space \space on $\partial\Omega_{\bar{\mathbf{f}}}$  \label{ContactForceNormalCond} \\
\bar{f}^{nor} \leq 0, & Non-tension, \qquad \space on $\partial\Omega_{\bar{\mathbf{f}}}$ \label{NonTensionNormalCond} \\
\textsl{g}^{nor} \geq 0, & Impenetrability, \space \space \space  on $\partial\Omega_{\bar{\mathbf{f}}}$ \label{ImpenetratibilityCond} \\
\textsl{g}^{nor} \bar{f}^{nor} = 0, & Complementarity, \space on $\partial\Omega_{\bar{\mathbf{f}}}$ \label{ComplementaryNormalCond}
\end{subnumcases}
whereas \EQ{CollinearityTangCond} to \EQ{ComplementaryTangCond} correspond to the tangential contact and friction laws, where the Coulomb friction model is adopted.
\begin{subnumcases}{}
\mathbf{e}^{tan}_{1} = - \mathbf{e}^{tan}_{2}, & Collinearity, \qquad \space \space on $\partial\Omega_{\bar{\mathbf{f}}}$ \label{CollinearityTangCond} \\
\bar{\mathbf{f}}_{1}^{tan} = - \bar{\mathbf{f}}_{2}^{tan}, & Collinearity, \qquad \space \space on $\partial\Omega_{\bar{\mathbf{f}}}$ \label{ContactForceTangCond} \\
| \bar{f}^{tan} | \leq \mu_{f} | \bar{f}^{nor} |, & Coulomb friction, \space on $\partial\Omega_{\bar{\mathbf{f}}}$ \label{CoulombCond} \\
| \textsl{g}^{tan} | \geq 0, & Slip/Non-Slip, \space \space \space \space \space \space  on $\partial\Omega_{\bar{\mathbf{f}}}$ \label{gsCond} \\
| \textsl{g}^{tan} | \left(  | \bar{f}^{tan} | - \mu_{f} | \bar{f}^{nor} |  \right) = 0 , & Complementarity, \space on $\partial\Omega_{\bar{\mathbf{f}}}$ \label{ComplementaryTangCond}
\end{subnumcases}

The Karush–Kuhn–Tucker (KKT) conditions \EQS{CollinearityNormalCond}-\eqref{ComplementaryTangCond} enforce normal and tangential contact laws at the boundary $\partial\Omega_{\bar{\mathbf{f}}}$. Symbols $\mathbf{e}^{nor}$ and $\mathbf{e}^{tan}$ denote the normal and tangential unit vectors, respectively. These constraints, illustrated in \FIG{fig:gov_eqns_frict_cont_probs:continuum_bodies_mpm_approximation}, uphold Newton's third law at the contact surface. The non-tension condition \EQ{NonTensionNormalCond} ensures non-stick contact.

\begin{figure}[htbp]
	\centering
	\begin{tabular}{ccc}
            \multicolumn{3}{c}{\includegraphics[scale=1.0]{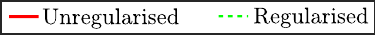}}\\
		\subfloat[\label{fig:gov_eqns_frict_cont_probs:contact_constraints:a}]{
			\includegraphics[scale=1.0]{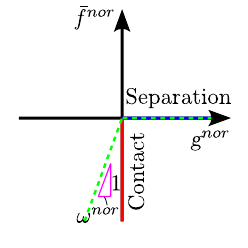}} &
            \subfloat[\label{fig:gov_eqns_frict_cont_probs:contact_constraints:b}]{
			\includegraphics[scale=1.0]{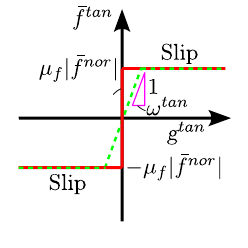}} &
            \subfloat[\label{fig:gov_eqns_frict_cont_probs:contact_constraints:c}]{
			\includegraphics[scale=1.0]{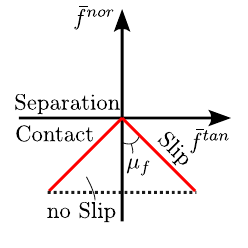}}
	\end{tabular}
	\caption[]{Kinematic contact constraints. \subref{fig:gov_eqns_frict_cont_probs:contact_constraints:a} normal contact law, \subref{fig:gov_eqns_frict_cont_probs:contact_constraints:b} tangential contact law and \subref{fig:gov_eqns_frict_cont_probs:contact_constraints:c} Coulomb’s cone for the two-dimensional problem. As a result of employing a penalty function approach and incorporating the penalty parameters $\omega^{nor}$ and $\omega^{tan}$, kinematic constraints are “weakly” imposed, as illustrated by a dashed green line.}
    \label{fig:gov_eqns_frict_cont_probs:contact_constraints}
\end{figure}

Additionally, impenetrability condition \EQ{ImpenetratibilityCond} prevents penetration between $\partial \Omega_{{1}\bar{\mathbf{f}}}$ and $\partial \Omega_{{2}\bar{\mathbf{f}}}$ when in contact. The Coulomb friction model in \EQ{CoulombCond} governs the tangential component $\bar{f}^{tan}$, related to slip contact. The condition $| \textsl{g}^{tan} | = 0$ denotes no slip, while $| \textsl{g}^{tan} | > 0$ indicates slip, satisfying $| \bar{f}^{tan} | = \mu_{f} | \bar{f}^{nor} |$. Regularisation through penalty parameters $\omega^{nor}$ and $\omega^{tan}$ weakly imposes these constraints. The green dashed line in \FIG{fig:gov_eqns_frict_cont_probs:contact_constraints:a} and \FIG{fig:gov_eqns_frict_cont_probs:contact_constraints:b} represents this weak enforcement, addressing potential zig-zagging in the tangential direction. For further details, refer to \cite{Liu2020ILS-MPM:Particles,DeLorenzis2011AAnalysis}. Section \ref{sec:definition_gap_function} elaborates on gap functions, slip functions, and contact force evaluation.

\section{An extended B-spline MPM for frictional contact} \label{sec:mpm_implementation}

\subsection{MPM approximation}

This work employs the Material Point Method (MPM) to discretise the governing equations \eqref{eqn:strong_form}. The domain $\Omega = \Omega_1 \cup \Omega_2$ is discretised into material points $\mathcal{Q} = { p | p = 1, 2,..., n_p }$. Let $n_{\mathcal{D}p}$ be the number of material points in discrete field $\mathcal{D}$ (e.g., $\Omega_{\mathcal{D}}$). For two solid bodies, $n_{p} = n_{1p} + n_{2p}$. \FIG{fig:gov_eqns_frict_cont_probs:continuum_bodies_mpm_approximation:b} illustrates the MPM for two separate discrete fields (master in red, slave in blue).

In MPM, mass density $\rho_{\mathcal{D}}$ and domain volume $\Omega_{\mathcal{D}}$ for discrete field $\mathcal{D}$ decompose into material point contributions per \EQS{eqn:mass_density_discrete_flds} and \eqref{eqn:volume_discrete_flds}:

\begin{equation}
\label{eqn:mass_density_discrete_flds}
\rho_{\mathcal{D}} \left(\mathbf{x}_{\mathcal{D}},t\right)=\sum\limits_{p = 1}^{{n_p}} {{\rho_{\mathcal{D}p}}\Omega_{\mathcal{D}p}\delta \left( {\mathbf{x}_{\mathcal{D}} - \mathbf{x}_{\mathcal{D}p}} \right)}
\end{equation}
and
\begin{equation}
\label{eqn:volume_discrete_flds}
\Omega_{\mathcal{D}} \left(\mathbf{x}_{\mathcal{D}},t\right)=\sum\limits_{p = 1}^{{n_p}} {{\Omega_{\mathcal{D}p}}\delta \left( {\mathbf{x}_{\mathcal{D}} - \mathbf{x}_{\mathcal{D}p}} \right)}, 
\end{equation}

Here, $\mathbf{x}_{\mathcal{D}}$ is the position vector of discrete field $\mathcal{D}$, and $\delta$ is the Dirac delta function. The material point mass density is $\rho_{\mathcal{D}p}=M_{\mathcal{D}p}/\Omega_{\mathcal{D}p}$, where $M_{\mathcal{D}p}$ and $\Omega_{\mathcal{D}p}$ and are the material point mass and volume, respectively. The symbol $\left( \cdot \right)_{\mathcal{D}p}$ denotes a quantity of material point $p$ at discrete field $\mathcal{D}$. Lagrangian material points move within a fixed Eulerian grid, consisting of $n_{n}$ grid nodes and $n_{cells}$ grid cells (\FIG{fig:gov_eqns_frict_cont_probs:continuum_bodies_mpm_approximation:b}).

\subsection{Discrete equilibrium equations}

The discrete equilibrium equations (\ref{eqn:strong_form}) are expressed for each discrete field $\mathcal{D}$ as
\begin{equation}
\begin{aligned}
\mathbf{M}_{\mathcal{D}} \ddot{\mathbf{u}}_{\mathcal{D}} + \mathbf{F}^{int}_{\mathcal{D}} = \mathbf{F}^{ext}_{\mathcal{D}} + \mathbf{F}^{cont}_{\mathcal{D}},
\end{aligned}
\label{eqn:EqOfMotion_Disc}
\end{equation}

with the last term representing externally applied contact forces. The lumped mass matrix $\mathbf{M}_{\mathcal{D}}$ is defined as
\begin{equation}
\begin{aligned}
M_{\mathcal{D}I} = \sum_{p=1}^{n_p} \Big( \rho_{\mathcal{D}p} N_{I}(\mathbf{x}_{p}) \Big) \Omega_{\mathcal{D}p}.
\end{aligned}
\label{eqn:MassLumpedMatrix}
\end{equation}

The nodal components of inertia forces $\mathbf{F}^{irt}_{\mathcal{D}I}$ and internal forces $\mathbf{F}^{int}_{\mathcal{D}I}$ are expressed as

$\mathbf{F}^{irt}_{\mathcal{D}I}$ are the nodal components of the inertia forces evaluated as
\begin{equation}
\label{eqn:InertiaForcesNodal}
\begin{aligned}
\mathbf{F}^{irt}_{\mathcal{D}I} = \sum_{p=1}^{n_p} (\rho_{\mathcal{D}p} \ddot{\mathbf{u}}_{\mathcal{D}p} \cdot N_I(\mathbf{x}_{p})) \Omega_{\mathcal{D}p}
\end{aligned}
\end{equation}
and
\begin{equation}
\label{eqn:InternalForcesNodal}
\begin{aligned}
\mathbf{F}^{int}_{\mathcal{D}I} = \sum_{p=1}^{n_p} (\boldsymbol{\sigma}_{\mathcal{D}p} \cdot \nabla N_I(\mathbf{x}_{p})) \Omega_{\mathcal{D}p}.
\end{aligned}
\end{equation}
, respectively.

The nodal components of the external force vector $\mathbf{F}^{ext}_{\mathcal{D}I}$ are computed as
\begin{equation}
\label{eqn:DiscreGalerkinDFGridArbi_3}
\begin{aligned}
\mathbf{F}^{ext}_I = \int_{\partial \Omega_{\bar{\mathbf{\tau}}}} ( \bar{\mathbf{\tau}} N_I(\mathbf{x}_{p})) \,d \partial \Omega_{\bar{\mathbf{\tau}}} + \sum_{p=1}^{n_p} \mathbf{b}_p N_I(\mathbf{x}_{p}) \Omega_p.
\end{aligned}
\end{equation}

The contact force nodal vector $\mathbf{F}^{cont}_{\mathcal{D}I}$ is defined as
\begin{equation}
\label{eqn:ContactForces}
\begin{aligned}
\mathbf{F}^{cont}_{\mathcal{D}I} = \int_{\partial \Omega_{\mathcal{D}\bar{\mathbf{f}}}} ( \bar{\mathbf{f}}_{\mathcal{D}}^{cont} N_I(\mathbf{x})) \,d \partial \Omega_{\mathcal{D}\bar{\mathbf{f}}}.
\end{aligned}
\end{equation}

The acceleration field is interpolated in a Galerkin sense as
\begin{equation}
\label{AccFunc}
\begin{aligned}
\ddot{\mathbf{u}}_{\mathcal{D}p} = \sum_{I=1}^{n_n} N_{I}(\mathbf{x}_{p}) \ddot{\mathbf{u}}_{\mathcal{D}I},
\end{aligned}
\end{equation}

where $\ddot{\mathbf{u}}_{\mathcal{D}I}$ represents the components of the nodal acceleration vector at node $I$. Displacement and velocity fields are derived similarly. \EQ{eqn:EqOfMotion_Disc} can be adapted to an explicit time integration scheme, as detailed in Section \ref{sec:solution_procedure}.

\subsection{Computational grid interpolation functions} \label{sec:comp_grid_interp_funcs}

In this study, the mapping between material points and grid nodes is achieved using Extended B-Spline interpolation functions, discussed in Section \ref{sec:eb_splines_interp_funcs}. These functions are derived from the standard B-Splines (referred to as Original B-Splines, or OBSs, in this work). Definition of the OBSs can be found in \cite{Hughes2005IsogeometricRefinement,Kakouris2018MaterialApproach}.

\subsubsection{Extended B-Spline interpolation functions} \label{sec:eb_splines_interp_funcs}

In this study, the approach from \cite{Yamaguchi2021ExtendedMethod} is adopted, leveraging OBS functions to minimise grid cell crossing errors, with EBS activated in boundary grid cells. EBS is crucial to mitigate stress errors arising from numerical integration issues in boundary cells. This is demonstrated in the numerical examples in Section \ref{sec:numerical_examples}. This section briefly discusses the EBS implementation to facilitate understanding of the classification of the basis functions as well as how EBS is fitted into frictional contact problems with MPM. The evaluation of Extended B-Splines involves a three-step procedure. \\
{\color{white}-}\\
\textbf{Step 1: Calculation of the volume fraction} \\
The initial step computes the grid cell volume fraction ($\phi_c$), representing the ratio of the total material point volume in a specific boundary grid cell to the cell's volume:
\begin{equation}
\begin{aligned}
\phi_c = \frac{\sum\limits_{p \in \mathcal{O}}^{} \Omega_p}{\Omega_c}.
\end{aligned}
\label{eqn:volume_fraction}
\end{equation}
{\color{white}-}\\
\textbf{Step 2: Classification of the grid cells} \\
Next, grid cells are classified as interior, boundary, or exterior based on the volume fraction:
\begin{equation}
\begin{aligned}
\text{Grid cell is categorised as}
    \begin{cases}
        \text{interior} & \text{if } \phi_c > C_c \\
        \text{boundary} & \text{if } 0 < \phi_c \leq C_c \\
        \text{exterior} & \text{if } \phi_c = 0
    \end{cases},
\end{aligned}
\label{eqn:grid_cell_int_bdy_ext}
\end{equation}
where $C_c$ is the \textit{occupation parameter}. \\
{\color{white}-}\\
\textbf{Step 3: Classification of basis functions} \\
The final step involves classifying B-Spline basis functions into stable, degenerate, and exterior groups:
\begin{itemize}
    \item \textbf{Stable:} if $\supp(N_I)$ contains at least one interior grid cell.
    \item \textbf{Degenerated:} if $\supp(N_I)$ contains no interior and at least one boundary grid cell.
    \item \textbf{Exterior:} if otherwise.
\end{itemize}

This classification ensures a robust representation of basis functions, addressing issues related to ill-conditioned matrices and non-uniform distribution of material points during large deformations. For more details, refer to \cite{Yamaguchi2021ExtendedMethod}.

\begin{figure}[htbp]
    \centering
    \includegraphics[scale=1.0]{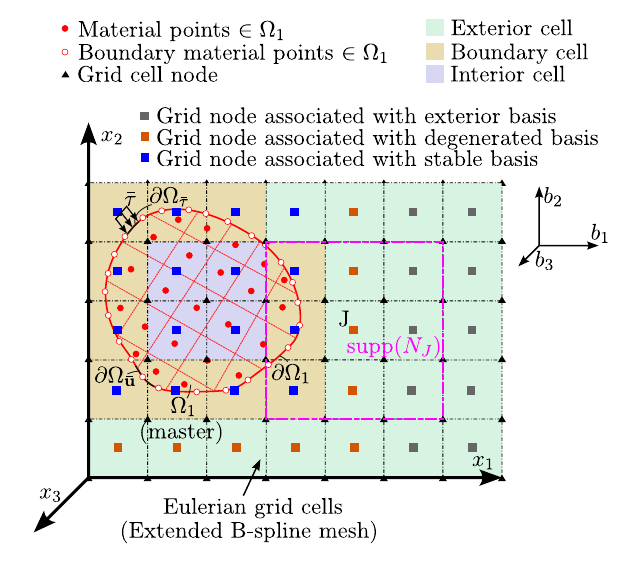}
	\caption{Categorising grid cells and nodes with classified bases on a two-dimensional continuum covered by an Eulerian grid. This is made according to the value of the volume fraction, $\phi_c$, at the boundary grid cells in relation to occupation parameter, $C_c$, chosen for the simulation.}
	\label{fig:mpm_implementation:classification_of_cell}
\end{figure}

In \FIG{fig:mpm_implementation:classification_of_cell}, a visual representation depicts the categorized bases, grid cells, and nodes using two-dimensional quadratic Extended B-Splines basis functions.

The stable and degenerated B-Splines, whose associated nodes are denoted by $\mathbb{I}$ and $\mathbb{J}$, respectively, are used to approximate the function near the boundary of a physical domain. The displacement field approximation is rewritten as
\begin{equation}
\label{DispFunc_ebs}
\begin{aligned}
\mathbf{u}_{\mathcal{D}p} = \sum_{I \in \mathbb{I}}^{} N_{I}(\mathbf{x}_{p}) \mathbf{u}_{\mathcal{D}I} + \sum_{J \in \mathbb{J}}^{} N_{J}(\mathbf{x}_{p}) \mathbf{u}_{\mathcal{D}J}
\end{aligned}
\end{equation}

\begin{figure}[htbp]
    \centering
    \includegraphics[scale=1.0]{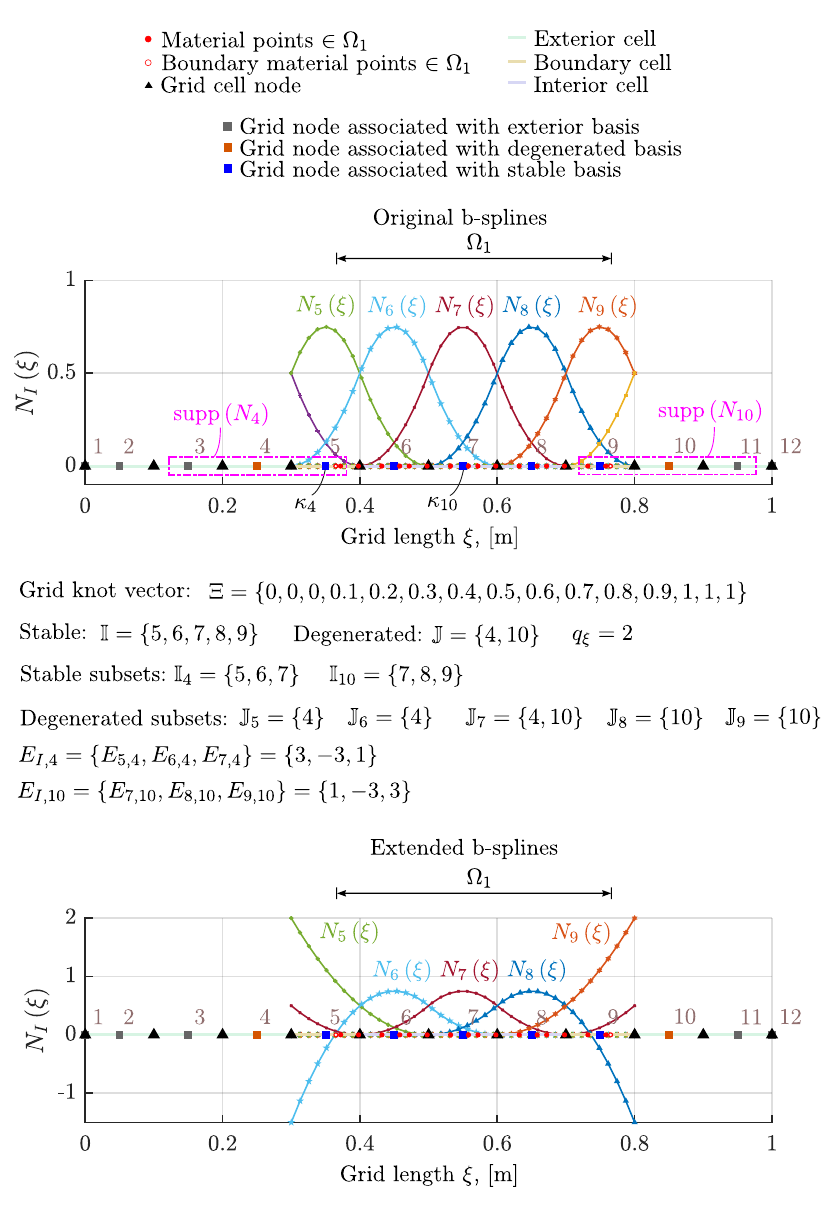}
	\caption{Derivation of quadratic Extended B-Splines from the Original B-Splines for the case of one-dimensional Eulerian grid. In this example, an occupation parameter at $C_c=0.75$ is utilised. Grid cells with identity, $4$ and $8$ are marked as boundary grid cells, resulting in two degenerated basis functions, associated with grid nodes (control points) $4$ and $10$.}
	\label{fig:mpm_implementation:1d_b_splines_ebs}
\end{figure}

Similar expressions are established for the velocity and acceleration fields. The coefficients of the degenerated bases are replaced by the linear combination of those of the stable bases as
\begin{equation}
\label{eqn:disp_J}
\begin{aligned}
\mathbf{u}_{J} = \sum_{I \in \mathbb{I}_{J}}^{} E_{IJ} \mathbf{u}_{I}
\end{aligned}
\end{equation}
where $E_{IJ}$ are the values of Lagrange polynomials at grid node $J$ for extrapolating the value of $\mathbf{u}$ at $J \in \mathbb{J}$ by using nodal values $\mathbf{u}_I$ with $I \in \mathbb{I}$. $\mathbb{I}_{J} \subseteq \mathbb{I}$ is a set of grid nodes associated with stable bases adjacent to grid node $J$. As already presented in \cite{Yamaguchi2021ExtendedMethod}, $E_{IJ}$ can be expressed as
\begin{equation}
\label{eqn:eij}
\begin{aligned}
E_{IJ} = \prod_{a=1}^{d} \left( \prod_{\chi=0,\chi \neq I^{a}-\kappa_{J}^{a}}^{q_{\xi}} \frac{J^a-\kappa_{J}^a-\chi}{I^a-\kappa_{J}^a-\chi} \right)
\end{aligned}
\end{equation}
when uniform knot vectors are used. The symbol, $\kappa$, denotes a starting point such as $\mathbb{I}_{J} = \kappa + [0,...,q_{\xi}]^{d}$. This is also illustrated in the one-dimensional example in \FIG{fig:mpm_implementation:1d_b_splines_ebs}.

The substitution of \EQ{eqn:eij} to \EQ{eqn:disp_J} and then to \EQ{DispFunc_ebs} yields to the following expression for interpolation
\begin{equation}
\label{eqn:disp_appox_ebs}
\begin{aligned}
\mathbf{u}_{\mathcal{D}p} = \sum_{I \in \mathbb{I}}^{} N^{ebs}_{I}(\mathbf{x}_{p}) \mathbf{u}_{\mathcal{D}I}
\end{aligned}
\end{equation}

The Extended B-Spline interpolation functions are eventually defined as
\begin{equation}
\label{eqn:shape_func_ebs}
\begin{aligned}
N^{ebs}_{I}(\mathbf{x}_{p}) = N^{}_{I}(\mathbf{x}_{p}) + \sum_{J \in \mathbb{J}_{I}}^{} N^{}_{J}(\mathbf{x}_{p}) E_{IJ},
\end{aligned}
\end{equation}
where $\mathbb{J}_{I}$ is the reciprocal set to $\mathbb{I}_{J}$. The first derivatives of the Extended B-Splines are expressed accordingly as
\begin{equation}
\label{eqn:der1_shape_func_ebs}
\begin{aligned}
\nabla N^{ebs}_{I}(\mathbf{x}_{p}) = \nabla N^{}_{I}(\mathbf{x}_{p}) + \sum_{J \in \mathbb{J}_{I}}^{} \nabla N^{}_{J}(\mathbf{x}_{p}) E_{IJ}.
\end{aligned}
\end{equation}

\FIG{fig:mpm_implementation:1d_b_splines_ebs} illustrate the derivation of quadratic Extended B-Splines from the Original B-Splines for the case of one-dimensional example. In this example, an occupation parameter, $C_c=0.75$ is utilised.

\subsection{Domain boundary representation} \label{sec:boundary_form}

In this work, solid body boundaries are represented by introducing additional material points, termed boundary material points. These points, depicted as ``unfilled" circles in figures e.g. \FIG{fig:gov_eqns_frict_cont_probs:continuum_bodies_mpm_approximation:b} and \FIG{fig:mpm_implementation:classification_of_cell}, serve to delineate the solid body boundaries. The total volume and mass of boundary material points is set to $0.1$\% of the domain's total volume and mass, respectively, in the numerical tests. Connecting these boundary material points with linear segments, as demonstrated in Section \ref{sec:numerical_examples}, yields good results. Although alternative boundary representations exist, such as those employing higher-order B-Splines \cite{Bing2019B-splineMethod}, they are beyond the scope of this work.

In the present work, the boundary material points are treated in a similar manner to the material points in the bulk; hence, no special treatment is required to update their properties, i.e., position, velocity, etc. The velocities of both the bulk and the boundary material points are projected onto the Eulerian grid to obtain the nodal values and solve the governing equations at grid nodes. The numerical implementation is detailed in Section \ref{sec:solution_procedure}.

In comparison, \cite{Guilkey2021AMethod} represented the boundary of the domain by defining a set of vertices - either line segments connected end to end in two dimensions or triangles in three dimensions. Regardless of dimensionality, these elements are massless and are advected using the velocity field of the associated material. However, special treatment is required for the massless elements at the domain boundary to follow bulk deformation when grid nodes with a non-zero value of the basis function have undefined velocity values. In particular, the velocities of the bulk material points are initially projected onto the Eulerian grid to obtain the nodal velocities. Next, the nodal velocities are updated and then are projected to the vertices that form the domain boundary. Vertices properties are updated utilising their velocity values. Although this process yields accurate estimates for frictional contact problems, it necessitates more computational time compared to the present approach, which updates the boundary material points in a more natural manner.

\subsection{Boundary tracking} \label{sec:boundary_track}

Boundary tracking, the process of identifying contact surfaces or points, is achieved in this work through a slave point - master segment pairing. The slave point is any boundary material point within the slave discrete field, while the master segment refers to a segment formed by consecutive boundary material points in the master discrete field.

\begin{figure}[htbp]
    \centering
    \includegraphics[scale=1.0]{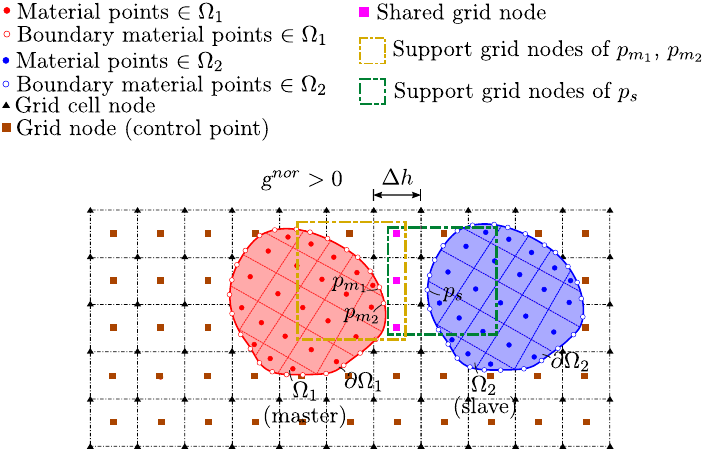}
	\caption{Boundary tracking process. Step 1: Shared control points indicate potential contact zones. Step 2: Distance of the boundary materials points are less than  grid spacing, $\Delta h$. Step 3: The $\textsl{g}^{nor} > 0$ which indicates that the slave material point $p_s$ is not in contact with the master segment, formed by boundary material points $p_{m_1}$ and $p_{m_2}$.}
	\label{fig:mpm_implementation:boundary_tracking}
\end{figure}

This contact pair detection involves a three-step process:
\begin{enumerate}
  \item \textbf{Common control points detection:} Check if boundary material points from different fields contribute to the same control point by projecting material point volumes to grid nodes using \EQ{eqn:nodal_volume}
  \begin{equation}
  \begin{aligned}
  \Omega_{\mathcal{D}I} = \sum_{p=1}^{n_p} \Big( N_{I}(\mathbf{x}_{p}) \Omega_{\mathcal{D}p} \Big).
  \end{aligned}
  \label{eqn:nodal_volume}
  \end{equation}   
  Shared control points indicate potential contact zones. For example, $\Omega_1$ and $\Omega_2$ share the same grid node $I$ when $\Omega_{1I}>0$ and $\Omega_{2I}>0$ simultaneously.
  \item \textbf{Distance check:} Examine if the distance between boundary material points from different fields is less than a minimum search size, typically the computational grid spacing denoted as $\Delta h$.
  \item \textbf{Contact pair determination:} Establish a slave boundary material point - master segment pair by verifying if $\textsl{g}^{nor} < 0$. This condition ensures the penetration of the slave boundary material point into the master segment. \FIG{fig:mpm_implementation:boundary_tracking} and \FIG{fig:mpm_implementation:contact_forces_2d} provide a visual representation of such pairs.
\end{enumerate}

The process described above results in great computational gains since detection of the slave-master segment pair is restricted to specific areas. This is a merit of the proposed method in compare to other MPM contact implementations, such as \cite{Guilkey2021AMethod}, which require massive intersection detection between the geometry components at each time step.

\subsection{Definition of the gap functions and contact forces evaluation} \label{sec:definition_gap_function}

To enforce contact constraints and prevent object penetration, a gap function measures the Euclidean distance between the boundaries of potentially contacting objects. This work adopts the normal gap and tangential slip functions introduced in \cite{Hamad2017AMPM}, which are expressed as.
\begin{equation}
\label{eqn:gap_function_normal}
\begin{aligned}
\prescript{(t)}{}{\textsl{g}}^{nor}_{} = \left( \prescript{(t)}{}{\mathbf{x}}^{}_{s} - \prescript{(t)}{}{\mathbf{x}}^{}_{\bar{s}} \right) \cdot \prescript{(t)}{}{\mathbf{e}}^{nor^T}_{} \quad \quad \text{on } \partial\Omega_{\bar{\mathbf{f}}}
\end{aligned}
\end{equation}
and
\begin{equation}
\label{eqn:gap_function_tangential}
\begin{aligned}
\prescript{(t)}{}{\textsl{g}}^{tan}_{} = \prescript{(t-\Delta t)}{}{l}^{}_{ms} \left( \prescript{(t)}{}{\beta}^{}_{\bar{s}} - \prescript{(t-\Delta t)}{}{\beta}^{}_{\bar{s}} \right)
\quad \quad \text{on } \partial\Omega_{\bar{\mathbf{f}}}
\end{aligned}
\end{equation}
, respectively. The symbol $l_{ms}$ denotes the length of the master segment while $\Delta t$ is the time step.

\begin{figure}[htbp]
	\centering
	\begin{tabular}{cc}
		\subfloat[\label{fig:mpm_implementation:contact_forces_2d:a}]{
			\includegraphics[scale=0.78]{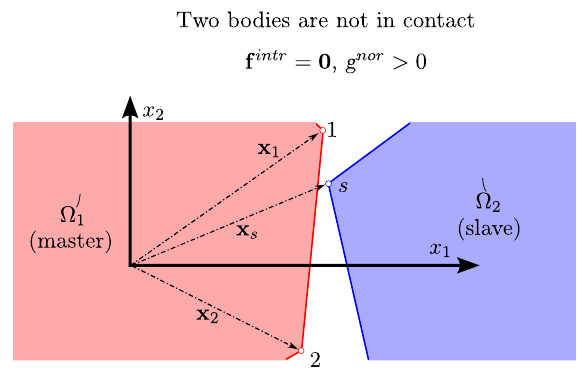}} &
            \subfloat[\label{fig:mpm_implementation:contact_forces_2d:b}]{
			\includegraphics[scale=0.78]{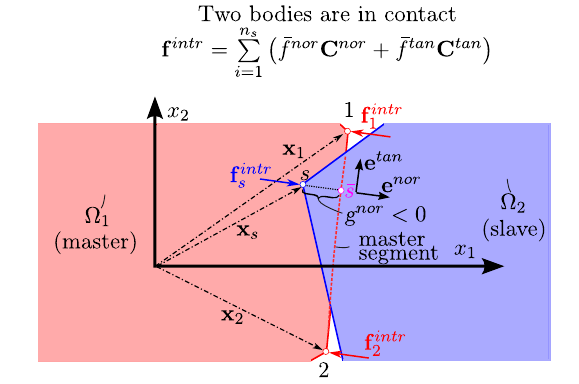}} \\
            \multicolumn{2}{c}{\subfloat[\label{fig:mpm_implementation:contact_forces_2d:c}]{
			\includegraphics[scale=0.78]{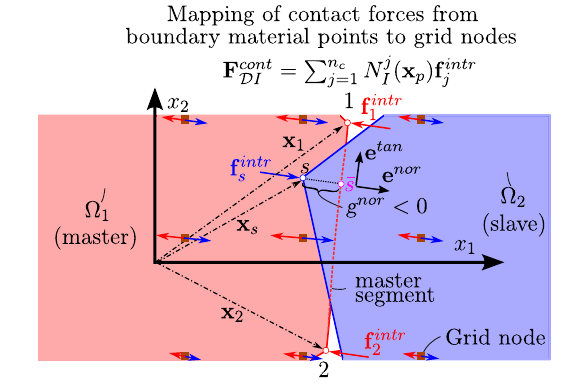}}}
	\end{tabular}
	\caption[]{\subref{fig:mpm_implementation:contact_forces_2d:a} Two bodies are not in contact. There are not interaction forces between the two bodies and normal gap function $\textsl{g}^{nor} >0$. \subref{fig:mpm_implementation:contact_forces_2d:b} two bodies are in contact. There are interaction forces between the two bodies. The slave boundary material point has penetrated the master segment. \subref{fig:mpm_implementation:contact_forces_2d:c} mapping of contact forces from material points to grid nodes. Then, the grid nodal contact forces are added to the grid nodal governing equations of motion (see \EQ{eqn:EqOfMotion_Disc}) which are updated.}
    \label{fig:mpm_implementation:contact_forces_2d}
\end{figure}

In \EQ{eqn:gap_function_normal}, $\mathbf{x}_{s}$ denotes the position vector of the slave node $s$, and $\mathbf{x}_{\bar{s}}$ is the position vector of $s$ projected onto the master segment. The normal gap function, $\textsl{g}^{nor}{}$, measures the Euclidean distance between the slave node $s$ and its projection $\bar{s}$ on the master segment. As discussed in Section \ref{sec:strong_form}, $\textsl{g}^{nor}{} > 0$ holds when the bodies are not in contact (see \FIG{fig:mpm_implementation:contact_forces_2d:a}), and $\textsl{g}^{nor}{} < 0$ when the slave body penetrates the master body (see \FIG{fig:mpm_implementation:contact_forces_2d:b} and \FIG{fig:mpm_implementation:contact_forces_2d:c}). In \FIG{fig:mpm_implementation:contact_forces_2d}, note that $\mathbf{e}^{nor}$ represents the unit outwards normal vector at the contact surface, starting from ${\mathbf{x}}^{}_{s}$ and ending at ${\mathbf{x}}^{}_{\bar{s}}$.

\begin{figure}[htbp]
	\centering
	\begin{tabular}{cc}
		\subfloat[\label{fig:mpm_implementation:tang_contact_forces_2d:a}]{
			\includegraphics[scale=0.83]{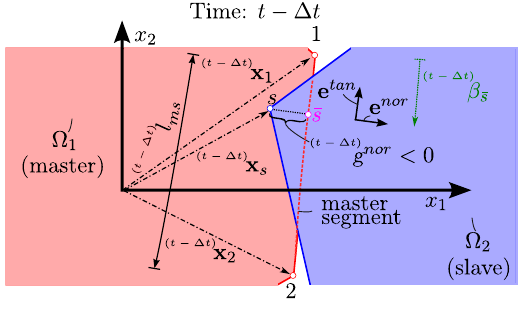}} &
            \subfloat[\label{fig:mpm_implementation:tang_contact_forces_2d:b}]{
			\includegraphics[scale=0.83]{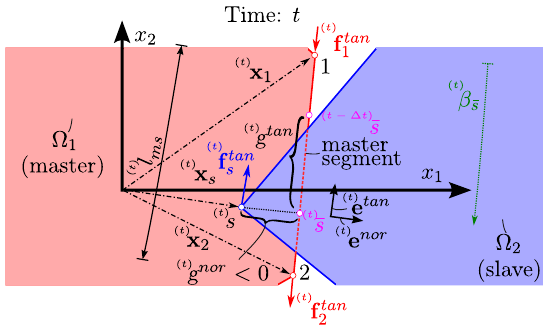}} \\
	\end{tabular}
	\caption[]{The tangential slip function measures the amount of relative movement of the two contact surfaces, $\partial \Omega_{1\bar{\mathbf{f}}}$ and $\partial \Omega_{2\bar{\mathbf{f}}}$, within a time increment. It is evaluation as per \EQ{eqn:gap_function_tangential}.The tangential slip produces the tangential contact forces at the slave boundary material points, $\mathbf{f}^{tan}_{s}$, and the master segment, $\mathbf{f}^{tan}_{1}$ and $\mathbf{f}^{tan}_{2}$.}
    \label{fig:mpm_implementation:tang_contact_forces_2d}
\end{figure}

In \EQ{eqn:gap_function_tangential} and for the tangential direction, $\mathbf{x}_1$ indicates the position vector of one end of the master segment. As also presented in \FIG{fig:mpm_implementation:contact_forces_2d}, position vectors $\mathbf{x}_1$ and $\mathbf{x}_2$ form the two ends of the master segment. The symbol, $\beta$, denotes the natural coordinate of $\bar{s}$ on the master segment and it is expressed as
\begin{equation}
\label{eqn:xi_natural_coord}
\begin{aligned}
\prescript{(t)}{}{\beta}^{}_{\bar{s}} = \frac{1}{\prescript{(t)}{}{l}^{}_{ms}} \left( \prescript{(t)}{}{\mathbf{x}}^{}_{\bar{s}} - \prescript{(t)}{}{\mathbf{x}}^{}_{1} \right)^T \cdot \prescript{(t)}{}{\mathbf{e}}^{tan}_{} \qquad \space \text{ and } \qquad \space 0<\beta<1
\end{aligned}
\end{equation}
The tangential slip function measures the amount of relative movement of the two contact surfaces, $\partial \Omega_{1\bar{\mathbf{f}}}$ and $\partial \Omega_{2\bar{\mathbf{f}}}$, within a time increment. A graphical representation of the tangential slip is illustrated in \FIG{fig:mpm_implementation:tang_contact_forces_2d} for two time steps, i.e. \subref{fig:mpm_implementation:tang_contact_forces_2d:a} $t - \Delta t$ and \subref{fig:mpm_implementation:tang_contact_forces_2d:b} $t$.

When the two bodies are in contact (see \FIG{fig:mpm_implementation:contact_forces_2d:b}), the normal contact component is evaluated from the \EQ{eqn:normal_contact_component} below
\begin{equation}
\prescript{(t)}{}{ \bar{f} }^{nor}_{} = \omega^{nor} \prescript{(t)}{}{ \textsl{g} }^{nor}_{}
\label{eqn:normal_contact_component}
\end{equation}
while the corresponding tangential component from \EQ{eqn:tangential_contact_component}
\begin{equation}
\prescript{(t)}{}{ \bar{f} }^{tan}_{} = \min \big( \mu_f | \bar{f}^{nor} | , \omega^{tan} | \prescript{(t)}{}{ \textsl{g}}^{tan}_{} | \big)  \text{sign} (-\prescript{(t)}{}{ \textsl{g}}^{tan}_{}).
\label{eqn:tangential_contact_component}
\end{equation}

\begin{Remark} In this work, an isotropic Coulomb friction model has been adopted to facilitate a direct comparison of the proposed method with other MPM variants in standard benchmark tests. Investigating more involved contact laws is beyond the scope of this research. We note however that an alternative friction contact law may be applied within the context of the proposed method by revisiting the definitions of the normal gap and tangential slip functions, and the tangential contact force component. The interested reader can refer to, e.g., \cite{Guilkey2023CohesiveMethod,Xiao2021DP-MPM:Fragmentation} for further insights.
\end{Remark}

Considering \EQS{eqn:normal_contact_component} and \eqref{eqn:tangential_contact_component}, the normal and tangential contact force vector are expressed as 
\begin{equation}
\prescript{(t)}{}{ \mathbf{f}}^{nor}_{} = \prescript{(t)}{}{ \bar{f} }^{nor}_{} \prescript{(t)}{}{\mathbf{C}}^{nor}_{}
\label{eqn:normal_contact_forces}
\end{equation}
and 
\begin{equation}
\prescript{(t)}{}{ \mathbf{f}}^{tan}_{} = \prescript{(t)}{}{ \bar{f} }^{tan}_{} \prescript{(t)}{}{\mathbf{C}}^{tan}_{},
\label{eqn:tangential_contact_forces}
\end{equation}
respectively, where the vectors ${\mathbf{C}}^{nor}$ and ${\mathbf{C}}^{tan}$ collect the action of the normal and tangential gap functions and are provided in Appendix \ref{App:Contact}. Eventually, the interaction force vector that emerges at the contact surface is expressed as
\begin{equation}
\prescript{(t)}{}{ \mathbf{f}}^{intr}_{} = \sum\limits_{i = 1}^{{n_s}} {\left(  \prescript{(t)}{}{ \mathbf{f}}^{nor}_{} + \prescript{(t)}{}{ \mathbf{f}}^{tan}_{} \right)}
\label{eqn:interaction_forces}
\end{equation}
where $n_s$ is the total number of slave nodes that have penetrated the master segment.

The interaction forces at the master and slave discrete field boundaries, introduced in \EQ{eqn:interaction_forces}, are then projected from the contact surface to the Eulerian grid according to the relation below
\begin{equation}
\prescript{(t)}{}{\mathbf{F}}^{cont}_{\mathcal{D}I} =\sum\limits_{j = 1}^{{n_c}} {N_I^j ( \prescript{(t)}{  }{ \mathbf{x} }^{  }_{ p } ) \prescript{(t)}{}{ \mathbf{f}}^{intr}_{j}}
\label{eqn:contact_forces_projection}
\end{equation}
where $n_c$ is the total number of contact grid nodes (see also \FIG{fig:mpm_implementation:contact_forces_2d:c}).

\begin{Remark}
We note that in the present work the boundary material points are not mass-less. In our view, this results in three interesting algorithim advantages compared to existing methods. First, the contact forces can be naturally projected directly from the slave - master segment pair onto the Eulerian grid at zero computational cost. Also, since the boundary material points have mass, the distribution of the contact forces to the computational nodes is not hindered, i.e., there is no requirement for redistributing the contact forces due to zero mass points. Finally, contrary to \cite{Guilkey2021AMethod}, the normal and tangential contact forces herein are explicitly evaluated from \EQS{eqn:normal_contact_forces} and \eqref{eqn:tangential_contact_forces}, respectively. Both contact components are projected from the slave-master segment pair to the computational grid. The governing equations are solved without any computational losses.
\end{Remark}

\subsection{Solution procedure} \label{sec:solution_procedure}

In this work, an explicit time integration scheme is utilised to integrate the equations of motion. Considering, a forward Euler integration scheme and a momentum formulation of the material point method as detailed in \cite{Kakouris2019Phase-FieldEnergy}, \EQ{eqn:EqOfMotion_Disc} is rewritten as
\begin{equation}
\begin{aligned}
\prescript{(t+\Delta t)}{}{\mathbf{p}}^{}_{\mathcal{D}I} = \prescript{(t)}{}{\mathbf{p}}^{}_{\mathcal{D}I} + \Delta t \left( \prescript{(t)}{}{\mathbf{F}}^{ext}_{\mathcal{D}I} + \prescript{(t)}{}{\mathbf{F}}^{cont}_{\mathcal{D}I} - \prescript{(t)}{}{\mathbf{F}}^{int}_{\mathcal{D}I} \right)
\end{aligned}
\label{eqn:ExpEqMotion_2}
\end{equation}
where $\prescript{(t+\Delta t)}{}{\mathbf{p}}^{}_{\mathcal{D}I}$ and $\prescript{(t)}{}{\mathbf{p}}^{}_{\mathcal{D}I}$ are the  nodal momentum at time $t+\Delta t$ and $t$, respectively. The nodal momentum at time $t$ is computed from the material points' mapping as 
\begin{equation}
\begin{aligned}
\prescript{(t)}{}{\mathbf{p}}^{}_{\mathcal{D}I} = \prescript{(t)}{}{\mathbf{M}}^{}_{\mathcal{D}I} \prescript{(t)}{}{\dot{\mathbf{u}}}^{}_{\mathcal{D}I} = \sum_{p=1}^{n_p} \Big( N_{I}(\prescript{(t)}{}{\mathbf{x}}^{}_{\mathcal{D}p}) \prescript{}{}{M}^{}_{\mathcal{D}p} \prescript{ (t)}{}{\dot{\mathbf{u}}}^{}_{\mathcal{D}p} \Big).
\end{aligned}
\label{eqn:momentum_t}
\end{equation}

The solution procedure of the proposed improved MPM contact algorithm is summarised in algorithm \ref{alg:mpm_implementation:ebs_mpm_frict_cont_probs}.  The first step is to define the initial state for the material points and the computational domain. As indicated in Section \ref{sec:definition_gap_function}, discrete fields, specifying subsets of material points, are labeled a priori and remain constant to enable tracking of their contact features during the analysis. Additionally, discrete field pairs (e.g., $\Omega_1$ - $\Omega_2$ pairs) are established. Prior to analysis, penalty parameters $\omega^{nor}$ and $\omega^{tan}$, along with the friction coefficient $\mu_f$, must be defined for each discrete field pair.

As indicated in Section \ref{sec:definition_gap_function}, discrete fields, specifying subsets of material points, are labeled a priori and remain constant to enable tracking of their contact features during the analysis. Additionally, discrete field pairs (e.g., $\Omega_1$ - $\Omega_2$ pairs) are established. Prior to analysis, penalty parameters $\omega^{nor}$ and $\omega^{tan}$, along with the friction coefficient $\mu_f$, must be defined for each discrete field pair.

In all numerical tests in this study, a structured computational grid is employed. For instance, a two-dimensional Eulerian grid is defined by the bottom-left coordinates $\mathbf{x}^{min} = \{ x^{min}_1, x^{min}_2 \}$ and the top-right coordinates $\mathbf{x}^{max} = \{ x^{max}_1, x^{max}_2 \}$. As emphasised in \cite{Nguyen2023TheApplications}, the use of a structured computational grid facilitates the identification of material points within their parent cells, ensuring robust computational efficiency. The grid spacing, denoted as $\Delta h$, is assumed to be uniform across all spatial dimensions. The time step, $\Delta t$, for the simulation is chosen following the methodology in Section \ref{sec:numerical_examples}.

After the initialisation step, \EQ{eqn:ExpEqMotion_2} is numerically solved until the desired simulation time. Time is discretised over $n_{steps}$ time increments. Furthermore, as detailed in algorithm \ref{alg:mpm_implementation:ebs_mpm_frict_cont_probs}, the basis functions are evaluated prior to grid mapping. Their evaluation process is presented in algorithm \ref{alg:mpm_implementation:ebs_basis_functions} and it is based on the three-step procedure introduced in Section \ref{sec:eb_splines_interp_funcs}. The equations of motion, \EQ{eqn:ExpEqMotion_2}, is then updated, where the contact forces are computed as outlined in Section \ref{sec:definition_gap_function}, \EQ{eqn:interaction_forces} and 
\EQ{eqn:contact_forces_projection}.

Next, the updated solution needs to be mapped back to material points. In this work, this is implemented utilising a so called Modified Update Stress Last (MUSL) approach \cite{Sulsky1995ApplicationMechanics, Nguyen2023TheApplications} to avoid numerical instabilities associated with potentially small grid nodal masses. It should be noted that the boundary material points are treated in the same manner as the material points in the bulk, e.g. position, strain, stress update etc. The procedure is presented in algorithm \ref{alg:mpm_implementation:map_from_grid_nodes_to_mps}.

\begin{Remark} 
In this work material points update uses a Particle-In-Cell approach. To this point, two main methods are used to update the velocity and position of material points in the literature, i.e. Particle In Cell (PIC) and Fluid Implicit Particle (FLIP) (see, e.g., \cite{Hammerquist2017AStability}). A PIC update filters velocity in each time step, which causes unwanted numerical diffusion, while FLIP eliminates that diffusion but may retain too much noise. Noise reduction via null space removal techniques attempts to minimise these errors \cite{Hammerquist2017AStability, Nairn2020NewFiltering}. In \cite{Jiang2017AnMethod} the affine PIC method has been introduced to reduce diffusion and improve conservation. Nevertheless, energy conservation and updating the velocity and position of material points remain open research issues. Very recently, \cite{Pretti2023AAnalysis} have proposed a novel time-stepping approach based on an efficient approximation of the Courant-Friedrich-Lewy (CFL) condition to enhance momenta and energy conservation in the MPM. We point however that such instabilities are different to the ones examined in this work and which pertain to the algorithmic treatment of frictional contact. However, revisiting the proposed formulation and further extending it using the aforementioned approaches is an interesting research direction that we would like to consider as future work.
\end{Remark}

\begin{algorithm}
    \caption{Extended B-Splines based MPM for frictional contact problems.}
    \begin{algorithmic}[1]
        \STATE \textbf{Data:} Define material point properties, computational grid, dynamic parameters ($\prescript{ (0)  }{   }{   \mathbf{x}     }^{   }_{  \mathcal{D}p }$, $\prescript{ (0)  }{   }{   \Omega     }^{   }_{  \mathcal{D}p }$, $\prescript{ (0)  }{   }{   \mathbf{u}     }^{   }_{  \mathcal{D}p }$, $\prescript{ (0)  }{   }{   \dot{\mathbf{u}}     }^{   }_{  \mathcal{D}p }$, $\prescript{ (0)  }{   }{   \ddot{\mathbf{u}}     }^{   }_{  \mathcal{D}p }$, $\prescript{ (0)  }{   }{   \boldsymbol {\varepsilon}     }^{   }_{  \mathcal{D}p }$, $\prescript{ (0)  }{   }{   \boldsymbol{\sigma}     }^{   }_{  \mathcal{D}p }$, $\prescript{ (0)  }{   }{   \rho     }^{   }_{  \mathcal{D}p }$, $E_{\mathcal{D}p}$, $\nu_{\mathcal{D}p}$, $\mathbf{x}^{min}$, $\mathbf{x}^{max}$, $\Delta h$, $t=0$, $\Delta t$).
        \STATE For each discrete field pair define: $\omega^{nor}$, $\omega^{tan}$ and $\mu_f$.
        \FOR{each time step $m=0,..,n_{steps}-1$ }
            \STATE  $t = t + \Delta t$.
            \STATE Reset grid quantities: $\prescript{(m)}{}{M}^{}_{\mathcal{D}I}=0$, $\prescript{(m)}{}{\mathbf{p}}^{}_{\mathcal{D}I} = \mathbf{0}$, $\prescript{(m)}{}{\mathbf{F}}^{int}_{\mathcal{D}I} = \mathbf{0}$,  $\prescript{(m)}{}{  \mathbf{F}}^{ext}_{ \mathcal{D}I} = \mathbf{0}$ and $\prescript{(m)}{}{\mathbf{F}}^{cont}_{\mathcal{D}I} = \mathbf{0}$, for each grid node $I$.
            \STATE Evaluate bases functions (see algorithm \ref{alg:mpm_implementation:ebs_basis_functions}).
            \STATE Map mass, momentum and internal forces from material points to grid nodes: $\prescript{(m)}{}{M}^{u}_{\mathcal{D}I}$, $\prescript{(m)}{}{\mathbf{p}}^{}_{\mathcal{D}I}$, $\prescript{(m)}{}{\mathbf{F}}^{int}_{\mathcal{D}I}$ and $\prescript{(m)}{}{\mathbf{F}}^{ext}_{\mathcal{D}I}$ (see \EQ{eqn:MassLumpedMatrix}, \EQ{eqn:momentum_t}, \EQ{eqn:InternalForcesNodal} and \EQ{eqn:DiscreGalerkinDFGridArbi_3})
            \STATE Apply (Dirichlet) boundary conditions: $\prescript{(m)}{}{p}^{}_{\mathcal{D}I,i } = 0$, $\prescript{(m)}{}{F}^{int}_{\mathcal{D}I,i} = 0$, $\prescript{(m)}{}{F}^{ext}_{\mathcal{D}I,i} = 0$ and $\prescript{(m)}{}{F}^{cont}_{\mathcal{D}I,i} = 0$, for each prescribed grid node $I$ in direction $i$.
            \STATE Update momentum: $\prescript{ (m+1) }{ }{ \widetilde{\mathbf{p}} }^{ }_{ \mathcal{D}I  }$ (see Eq. \ref{eqn:ExpEqMotion_2}) \footnotemark{} 
            \STATE Map from grid nodes to material points (see algorithm \ref{alg:mpm_implementation:map_from_grid_nodes_to_mps}).  \STATE Export output (post-processing).  
        \ENDFOR
    \end{algorithmic}
	\label{alg:mpm_implementation:ebs_mpm_frict_cont_probs}
\end{algorithm}
\footnotetext{The symbol $\sim$ indicates trial (temporal) quantities.}

\begin{algorithm}
    \caption{Evaluation of Extended B-Spline basis functions.}
    \begin{algorithmic}[1]
        \STATE Calculate grid cells' volume fraction (see \EQ{eqn:volume_fraction}).
        \STATE Classify grid cells into interior, boundary, and exterior ones according to their material point volume fractions as compared to occupation parameter (see \EQ{eqn:grid_cell_int_bdy_ext}).
        \STATE Categorise bases functions as either stable, degenerated, or exterior (see step three in Section \ref{sec:eb_splines_interp_funcs}).
        \STATE Determine subsets $\mathbb{I}_J$ for degenerated bases functions.       
        \STATE Compute: $\mathbf{N}_{} (\prescript{(m)}{}{ \mathbf{x} }^{}_{p} )$, $\nabla \mathbf{N}_{} ( \prescript{(m)}{}{ \mathbf{x} }^{}_{p} )$, for all material points. These correspond to the Original B-Splines as presented in \cite{Hughes2005IsogeometricRefinement,Kakouris2018MaterialApproach}.
        \STATE Compute: $\mathbf{N}_{}^{ebs} (\prescript{(m)}{}{ \mathbf{x} }^{}_{p} )$, $\nabla \mathbf{N}_{}^{ebs} ( \prescript{(m)}{}{ \mathbf{x} }^{}_{p} )$, only for degenerated basis functions (see \EQ{eqn:shape_func_ebs} and \EQ{eqn:der1_shape_func_ebs}).
    \end{algorithmic}
	\label{alg:mpm_implementation:ebs_basis_functions}
\end{algorithm}

\begin{algorithm}
    \caption{Map from grid nodes to material points. It is based on Modified Update Stress Last (MUSL) algorithm \cite{Sulsky1995ApplicationMechanics, Nguyen2023TheApplications}.}
    \begin{algorithmic}[1]
        \STATE \textbf{Update material points' position, displacement, velocity and grid node velocities (double mapping)}        
        \bindent
        \STATE Get nodal velocities and accelerations: $\prescript{(m+1)}{}{ \dot{\widetilde{\mathbf{u}}}}^{}_{\mathcal{D}I} = \prescript{(m+1)}{}{ \widetilde{\mathbf{p}}}^{ }_{\mathcal{D}I}/\prescript{(m)}{}{M}^{}_{\mathcal{D}I}$ and $\prescript{(m+1)}{}{\ddot{\widetilde{\mathbf{u}}}}^{}_{\mathcal{D}I} = (\prescript{(m)}{}{{\mathbf{F}} }^{ext}_{\mathcal{D}I} + \prescript{(m)}{}{{\mathbf{F}} }^{cont}_{\mathcal{D}I} - \prescript{(m)}{}{{\mathbf{F}} }^{int}_{\mathcal{D}I}) /\prescript{(m)}{}{M}^{}_{\mathcal{D}I}$.
        \STATE \begin{flushleft} Update material point position: $\prescript{ (m+1)}{}{\mathbf{x}}^{}_{\mathcal{D}p} = \prescript{(m)}{}{\mathbf{x} }^{}_{\mathcal{D}p} + \Delta t \sum_{I=1}^{n_n} N_{I}(\prescript{(m)}{}{ \mathbf{x}}^{}_{\mathcal{D}p}) \prescript{(m+1)}{}{\dot{\widetilde{\mathbf{u}}}}^{}_{\mathcal{D}I}$. \end{flushleft}
        \STATE \begin{flushleft} Update material point velocity: $\prescript{ (m+1)}{}{\dot{\mathbf{u}}}^{}_{\mathcal{D}p} = \sum_{I=1}^{n_n} N_{I}( \prescript{(m)}{}{\mathbf{x}}^{}_{\mathcal{D}p}) \prescript{(m+1)}{}{ \dot{\widetilde{\mathbf{u}}}}^{}_{\mathcal{D}I}$. \end{flushleft}
        \STATE \begin{flushleft} Update material point displacement: $\prescript{(m+1)}{}{\mathbf{u}}^{}_{\mathcal{D}p} = \prescript{(m+1)}{}{\mathbf{x} }^{}_{\mathcal{D}p} - \prescript{(0)}{}{\mathbf{x}}^{}_{\mathcal{D}p}$.
        \end{flushleft}
        \STATE \begin{flushleft} Map momentum from material points to grid nodes: $\prescript{(m+1)}{}{ \mathbf{p}}^{}_{\mathcal{D}I} = \sum_{p=1}^{n_p} N_{I}(\prescript{(m)}{}{ \mathbf{x}}^{}_{\mathcal{D}p}) \left( M_{\mathcal{D}p} \prescript{(m+1)}{}{ \dot{\mathbf{u}}}^{}_{\mathcal{D}p} \right)$.
        \end{flushleft}
        \STATE Apply (Dirichlet) boundary conditions: $\prescript{ (m+1) }{ }{ p }^{ }_{ \mathcal{D}I,i } = 0$, for each prescribed grid node $I$ in direction $i$.
        \eindent
        \STATE \textbf{end}
        \STATE \textbf{Update material points. From grid nodes to material points.}
        \bindent
        \STATE Get nodal velocities: $\prescript{ (m+1)  }{ }{ \dot{\mathbf{u}} }^{ }_{ \mathcal{D}I } = \prescript{ (m+1) }{ }{ \mathbf{p} }^{ }_{ \mathcal{D}I  } / \prescript{ (m)  }{ }{   M }^{ }_{ \mathcal{D}I  }$.
        \STATE \begin{flushleft} Compute material point gradient velocity: $\prescript{(m+1)}{}{ \mathbf{L}}^{}_{\mathcal{D}p} = \sum\limits_{I = 1}^{{n_n}} \nabla N_I \left(\prescript{ (m) }{ }{\mathbf{x}}^{ }_{\mathcal{D}p} \right) \prescript{ (m+1)  }{ }{ \dot{\mathbf{u}} }^{ }_{ \mathcal{D}I } $ \end{flushleft}
        \STATE \begin{flushleft} Update gradient deformation tensor: $\prescript{(m+1)}{}{ \mathbf{F}}^{}_{\mathcal{D}p} = \left(  \mathbf{I} + \prescript{(m+1)}{}{ \mathbf{L}}^{}_{\mathcal{D}p} \Delta t \right) \prescript{(m)}{}{ \mathbf{F}}^{}_{\mathcal{D}p}$ \end{flushleft}        
        \STATE \begin{flushleft} Get material point incremental strain: $\prescript{}{}{\boldsymbol {\Delta \varepsilon}}^{}_{\mathcal{D}p } = 0.5 \Delta t \left( \prescript{ (m+1)}{}{\mathbf{L}}^{ }_{\mathcal{D}p} + \prescript{(m+1)}{}{ \mathbf{L}}^{T}_{ \mathcal{D}p} \right)$
        \end{flushleft}        
        \STATE Update material point strains: $\prescript{ (m+1)  }{   }{   \boldsymbol {\varepsilon}     }^{   }_{  \mathcal{D}p } = \prescript{ (m)  }{   }{   \boldsymbol {\varepsilon}     }^{   }_{  \mathcal{D}p } + \prescript{  }{   }{   \boldsymbol {\Delta \varepsilon}     }^{   }_{  \mathcal{D}p }$
        \STATE Update material point stresses: $\prescript{ (m+1)  }{   }{   \boldsymbol{\sigma}     }^{   }_{  \mathcal{D}p } = \partial 
        \prescript{(m+1)}{}{ \psi    }^{   }_{  el\mathcal{D}p } / \partial          \prescript{ (m+1)  }{   }{   \boldsymbol {\varepsilon}     }^{   }_{  \mathcal{D}p }$.
        \eindent
        \STATE \textbf{end}
    \end{algorithmic}
	\label{alg:mpm_implementation:map_from_grid_nodes_to_mps}
\end{algorithm}

\section{Numerical examples} \label{sec:numerical_examples}

This section presents numerical validations of the proposed model, focusing on one and two-dimensional contact problems. Results are compared with literature and analytical solutions, including tests with OBS and state-of-the-art MPM contact implementations. A comparison between EBS and OBS is vital in frictional contact problems to illustrate the influence of numerical integration errors at the boundary grid cells and ultimately at the evolved contact forces. This is an aspect that has never been explored in frictional contact problems with other MPM variants.

Quadratic B-Splines ($C^{1}$) are employed for the background grid, ensuring an initial cell density of at least $4$ material points per grid cell in one-dimensional examples and $16$ material points per grid cell in two-dimensional problems. In all two-dimensional problems, plane strain conditions are assumed. In addition, linear elastic material response is considered for all benchmarks except from Section \ref{sec:imp_elast_rings} where a hyperelastic compressible Neo-Hookean material law has been utilised.

The explicit time integration scheme's stability is maintained by limiting the time increment, $\Delta t$, with the upper bound:
\begin{equation}
\Delta t \leq 0.1 \Delta t_{cr}
\label{eqn:time_step_condition}
\end{equation}
Here, $\Delta t_{cr} = \Delta h / c_0$ is the critical time step defined by the Courant–Friedrichs–Lewy (CFL) condition, where $c_0 = \sqrt{E/\rho}$ represents the speed of sound in solids.

To prevent contact violation, penalty parameters $\omega^{nor}$ and $\omega^{tan}$ need substantial values, yet excessively large ones can induce numerical instabilities. Thus, the following relation is employed to determine these parameters:
\begin{equation}
\omega^{nor} = \omega^{tan} = \frac{1}{\Delta h / E}
\label{eqn:penalty_parameter}
\end{equation}
This equation, used in previous works \cite{Kikuchi1988ContactMethods, Leichner2019ARepresentation, Liu2020ILS-MPM:Particles}, demonstrates stability in the conducted numerical tests, avoiding associated instabilities.

\subsection{One-dimensional compression of two contacting bars under self weight}
\label{sec:1d_comp2bars_under_self_weight}

\begin{figure}[htbp]
    \centering
    \includegraphics[scale=1.0]{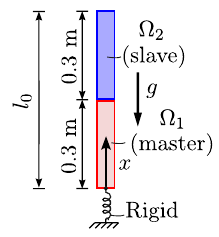}
	\caption{One-dimensional compression of two contacting bars under self weight. Geometry and boundary conditions.}
	\label{fig:1d_comp2bars_under_self_weight:geo_ibcs}
\end{figure}

This example investigates the one-dimensional compression of two bars under self-weight, as illustrated in \FIG{fig:1d_comp2bars_under_self_weight:geo_ibcs}. This serves to validate the proposed algorithm against analytical solutions and highlight its advantages over OBS.

The setup involves two bars, each of length $0.3$ m, positioned one on top of the other (see \FIG{fig:1d_comp2bars_under_self_weight:geo_ibcs}). The bottom bar is fixed at its lower edge with $u \left( x = 0 \right) = 0$. To impose this fixed boundary condition, a “rigid” spring with a stiffness of $65 \cdot 10^9$ N/m is added to the boundary material point. This ensures fixed boundary conditions even when the control point of the Eulerian grid does not coincide with the boundary material point. Both bars share identical material properties: mass density $\rho = 2783$ kg/m\textsuperscript{3}, Young's modulus $E = 50.5$ GPa, and cross-sectional area $A = 1$ m\textsuperscript{2}.

In this example, a one-dimensional Eulerian grid is defined with a grid spacing of $\Delta h = 0.1$ m. Each discrete field has initially a grid cell density of $4$ material points per cell, along with $2$ boundary material points. The time step, $\Delta t$, follows \EQ{eqn:time_step_condition}, set as $\Delta t = 0.1 \Delta t_{cr} = 2.34753 \cdot 10^{-6}$ sec. The bars are initially stress-free at $t=0$, and a quasi-static response is ensured by gradually increasing the gravity force to $g = -9.81$ m/s\textsuperscript{2} over $0.0281$ sec (equivalent to $12000$ time steps), $200$ times the elastic wave transit time in the bars, i.e., $(1/c_0) \cdot l_0 \cdot 200$. The EBS adopts an occupation parameter of $0.75$. The normal penalty parameter, calculated using \EQ{eqn:penalty_parameter}, is set to $\omega^{nor} = 505 \cdot 10^9$ N/m\textsuperscript{3}.

To examine the influence of the Eulerian grid on the evaluated stress as a result of incomplete integration, $6$ cases are considered as initial bar locations  (termed as bar offsets herein), i.e. $0$, $+0.016$, $+0.033$, $+0.05$, $+0.066$ and $+0.083$ m. The master bar is initially located at position $0.2$ m on the Eulerian grid and it is offset as per the bar offset cases.

The analytical solution for the Cauchy stress over bars' length, $\sigma$, are computed from the \EQ{eqn:1d_comp2bars_under_self_weight:anal_solution} below \cite{Zhang2011MaterialFunction}.
\begin{equation}
\begin{aligned}
\sigma \left( x \right) = \sigma \left( 0 \right) \frac{(1 - x/l_0) + \alpha (x/l_0)}{1 - \alpha (1-\alpha) (x/l_0)}
\end{aligned}
\label{eqn:1d_comp2bars_under_self_weight:anal_solution}
\end{equation}
where $\sigma \left( 0 \right) = \rho g A l_0$ and $\alpha = \sigma \left( 0 \right) / (2 E)$. Total weight of the bar is $W = -16380.738$ N.

\begin{figure}[htbp]
	\centering
	\begin{tabular}{ccc}
            \multicolumn{3}{c}{\includegraphics[scale=1.0]{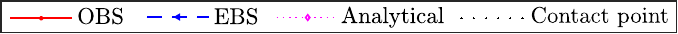}}\\
		\subfloat[Bar Offset: 0.000 m \label{fig:1d_comp2bars_under_self_weight:stress:a}]{
			\includegraphics[scale=1.0]{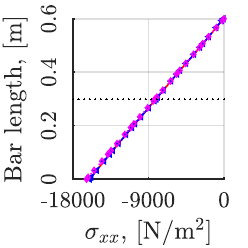}} &
            \subfloat[Bar Offset: 0.016 m \label{fig:1d_comp2bars_under_self_weight:stress:b}]{
			\includegraphics[scale=1.0]{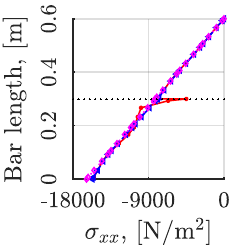}} & 
            \subfloat[Bar Offset: 0.033 m \label{fig:1d_comp2bars_under_self_weight:stress:c}]{
			\includegraphics[scale=1.0]{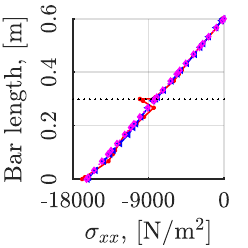}} \\
            \subfloat[Bar Offset: 0.050 m \label{fig:1d_comp2bars_under_self_weight:stress:d}]{
			\includegraphics[scale=1.0]{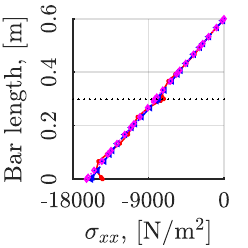}} &
            \subfloat[Bar Offset: 0.066 m \label{fig:1d_comp2bars_under_self_weight:stress:e}]{
			\includegraphics[scale=1.0]{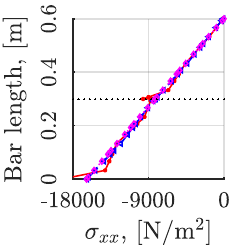}} & 
            \subfloat[Bar Offset: 0.083 m \label{fig:1d_comp2bars_under_self_weight:stress:f}]{
			\includegraphics[scale=1.0]{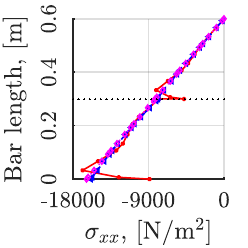}}
	\end{tabular}
	\caption[]{One-dimensional compression of two contacting bars under self weight. Stress over bar length for various bar offsets. The proposed frictional contact EBS-based algorithm (blue dashed line) provides good estimates of the bars' stress along their length, regardless the position of the bar with respect to the Eulerian grid. In comparison, OBS (red line) fails to converge to the accurate stress, showing severe errors at the location of the contact and bottom points.}
	\label{fig:1d_comp2bars_under_self_weight:stress}
\end{figure}

The stress computations, illustrated in \FIG{fig:1d_comp2bars_under_self_weight:stress} for both OBS and EBS against the analytical solution, highlight the accuracy of EBS in predicting stress across various bar offsets. Regardless of the bar's position relative to the Eulerian grid, EBS provides precise stress predictions. In comparison, OBS fails to converge accurately and exhibits substantial errors at the contact point. Stress noise and numerical integration errors are evident with OBS, especially at the master bar's bottom edge for certain bar offsets ($+0.050$, $+0.066$, and $+0.083$, see \FIG{fig:1d_comp2bars_under_self_weight:stress}).

The dimensionless error, calculated using \EQ{eqn:1d_comp_under_self_weight:error}, further emphasises these differences.
\begin{equation}
\text{Error} = \sum\limits_{p = 1}^{{n_p}} {\frac{ | \sigma^{\text{numerical}}_p - \sigma^{\text{analytical}}_p | \Omega_p}{|W| l_0}}
\label{eqn:1d_comp_under_self_weight:error}
\end{equation}

\begin{figure}[htbp]
    \centering
    \includegraphics[scale=1.0]{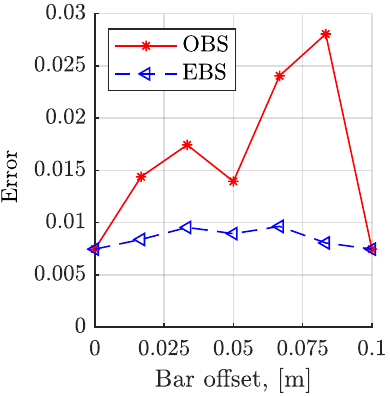}
	\caption{One-dimensional compression of two contacting bars under self weight. Stress error for various bar offsets within the Eulerian grid. The stress error is retained at low levels with the proposed EBS (blue dashed line). This is not the case for OBS which stress error increases with the bar offset due to incomplete numerical integration (red line).}
	\label{fig:1d_comp2bars_under_self_weight:error}
\end{figure}

As shown in \FIG{fig:1d_comp2bars_under_self_weight:error}, OBS error increases due to incomplete numerical integration, while EBS error remains consistently low.

\begin{figure}[htbp]
	\centering
	\begin{tabular}{cc}   
   \multicolumn{2}{c}{\includegraphics[scale=1.0]{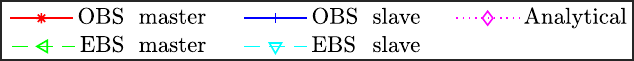}}\\
		\subfloat[\label{fig:1d_comp2bars_under_self_weight:bottom_middle_stress:a}]{
			\includegraphics[scale=1.0]{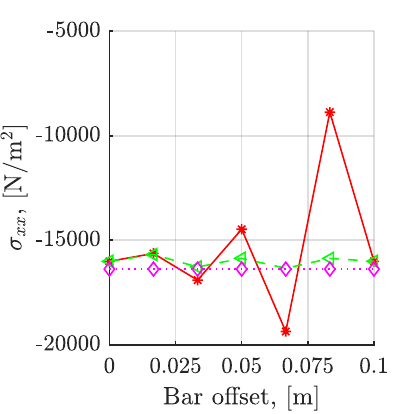}} &
            \subfloat[\label{fig:1d_comp2bars_under_self_weight:bottom_middle_stress:b}]{
			\includegraphics[scale=1.0]{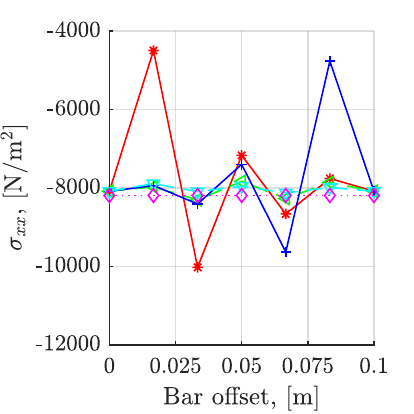}}
	\end{tabular}
	\caption[]{One-dimensional compression of two contacting bars under self weight. \subref{fig:1d_comp2bars_under_self_weight:bottom_middle_stress:a} Master bar bottom stress and  \subref{fig:1d_comp2bars_under_self_weight:bottom_middle_stress:b} Master bar top, slave bar bottom Master bar top stress, slave bar bottom. In all cases, the OBS results significantly fluctuate from the analytical solution while minor divergence is observed using the EBS.}
	\label{fig:1d_comp2bars_under_self_weight:bottom_middle_stress}
\end{figure}

Additionally, \FIG{fig:1d_comp2bars_under_self_weight:bottom_middle_stress:a} and \FIG{fig:1d_comp2bars_under_self_weight:bottom_middle_stress:b} present computed stresses at critical bar points for various bar offsets. \FIG{fig:1d_comp2bars_under_self_weight:bottom_middle_stress:a} reveals that OBS results fluctuate significantly from the analytical solution, while EBS shows minor divergence. Furthermore, \FIG{fig:1d_comp2bars_under_self_weight:bottom_middle_stress:b} illustrates contact stresses, where OBS not only deviates significantly from the analytical solution but also exhibits substantial differences between master top and slave bottom stresses, exceeding $50$\%. In contrast, EBS contact stresses show differences of up to $3$\%.

\subsection{Longitudinal impact of two bars}
\label{sec:longitudinal_impact_of_two_bars}

\begin{figure}[htbp]
    \centering
    \includegraphics[scale=1.0]{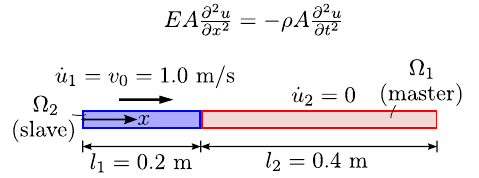}
	\caption{Longitudinal impact of two bars. Geometry and initial conditions.}
	\label{fig:longitudinal_impact_of_two_bars:geo_ibcs}
\end{figure}

The second example concerns the investigation of the longitudinal impact of two one-dimensional bars, characterised by different initial velocities, thereby resulting in a collision. The aim of this example is to ascertain the post-impact velocities of the objects, as well as the corresponding impact force over time. In the present scenario, a one-dimensional case is examined where a bar with length of $l_1 = 0.2$ m collides with another bar with length $l_2 = 0.4$ m as depicted in \FIG{fig:longitudinal_impact_of_two_bars:geo_ibcs}.

\begin{figure}[htbp]
	\centering
	\begin{tabular}{cc}
            \multicolumn{2}{c}{\includegraphics[scale=0.85]{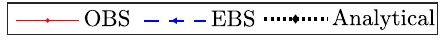}}\\
		  \subfloat[\label{fig:longitudinal_impact_of_two_bars:vel:a}]{
			\includegraphics[scale=0.85]{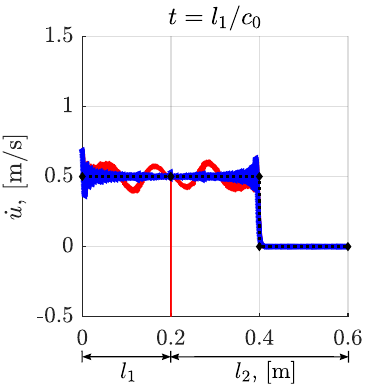}} &
            \subfloat[\label{fig:longitudinal_impact_of_two_bars:vel:b}]{
			\includegraphics[scale=0.85]{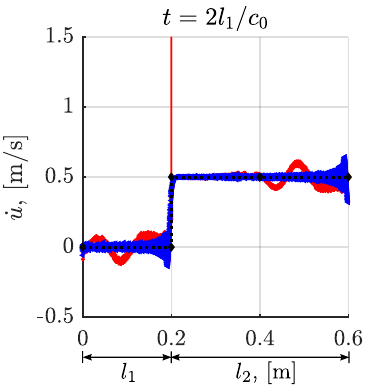}} \\ 
            \subfloat[\label{fig:longitudinal_impact_of_two_bars:vel:c}]{
			\includegraphics[scale=0.85]{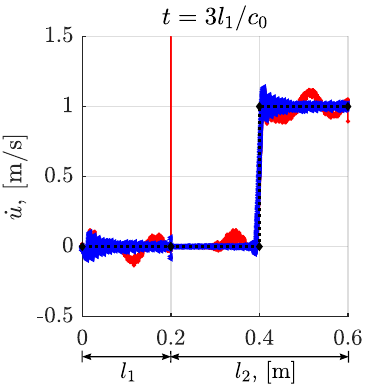}} &
            \subfloat[\label{fig:longitudinal_impact_of_two_bars:vel:d}]{
			\includegraphics[scale=0.85]{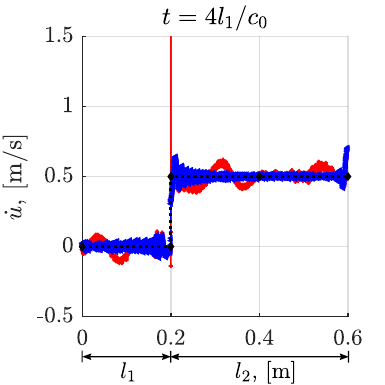}} \\
            \multicolumn{2}{c}{\subfloat[ \label{fig:longitudinal_impact_of_two_bars:vel:e}]{
			\includegraphics[scale=0.85]{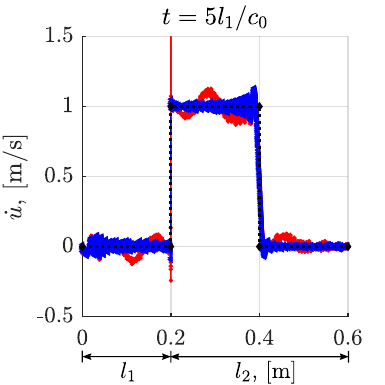}}} 
	\end{tabular}
	\caption[]{Longitudinal impact of two bars. Velocity over bars' length for time instances, \subref{fig:longitudinal_impact_of_two_bars:vel:a} $l_1/c_0 = 46.95$ $\mu$s, \subref{fig:longitudinal_impact_of_two_bars:vel:b} $2 l_1/c_0 = 93.9$ $\mu$s, \subref{fig:longitudinal_impact_of_two_bars:vel:c} $3 l_1/c_0 = 140.85$ $\mu$s,\subref{fig:longitudinal_impact_of_two_bars:vel:d} $4 l_1/c_0 = 187.8$ $\mu$s \subref{fig:longitudinal_impact_of_two_bars:vel:e} $5 l_1/c_0 = 187.8$ $\mu$s. Two numerical simulations are considered with OBS (red line) and EBS (blue dashed line). The analytical solutions are also plotted with black dotted line. EBS yields better accuracy over OBS and converges to the analytical solutions.}
	\label{fig:longitudinal_impact_of_two_bars:vel}
\end{figure}

\begin{figure}[htbp]
	\centering
	\begin{tabular}{cc}
            \multicolumn{2}{c}{\includegraphics[scale=0.85]{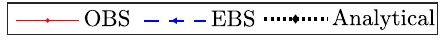}}\\
		  \subfloat[\label{fig:longitudinal_impact_of_two_bars:stress:a}]{
			\includegraphics[scale=0.85]{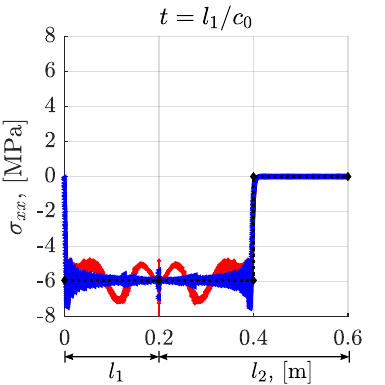}} &
            \subfloat[\label{fig:longitudinal_impact_of_two_bars:stress:b}]{
			\includegraphics[scale=0.85]{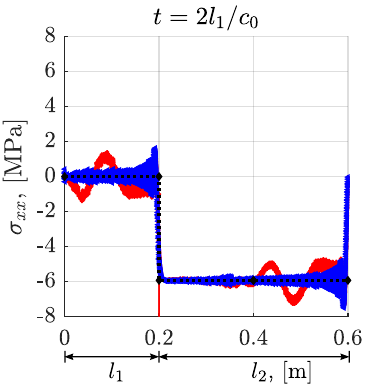}} \\ 
            \subfloat[\label{fig:longitudinal_impact_of_two_bars:stress:c}]{
			\includegraphics[scale=0.85]{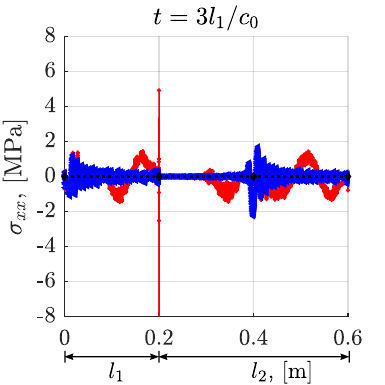}} &
            \subfloat[\label{fig:longitudinal_impact_of_two_bars:stress:d}]{
			\includegraphics[scale=0.85]{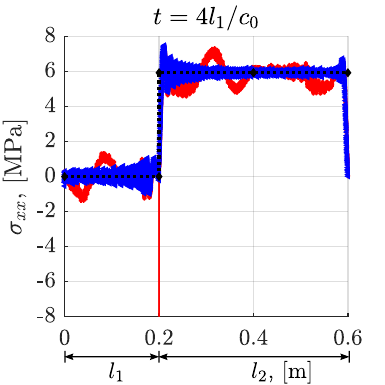}} \\
            \multicolumn{2}{c}{\subfloat[ \label{fig:longitudinal_impact_of_two_bars:stress:e}]{
			\includegraphics[scale=0.85]{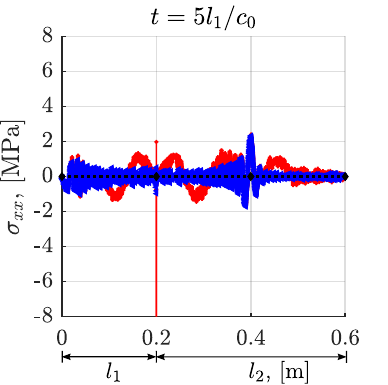}}} 
	\end{tabular}
	\caption[]{Longitudinal impact of two bars. Stress over bars' length for time instances, \subref{fig:longitudinal_impact_of_two_bars:stress:a} $l_1/c_0 = 46.95$ $\mu$s, \subref{fig:longitudinal_impact_of_two_bars:stress:b} $2 l_1/c_0 = 93.9$ $\mu$s, \subref{fig:longitudinal_impact_of_two_bars:stress:c} $3 l_1/c_0 = 140.85$ $\mu$s,\subref{fig:longitudinal_impact_of_two_bars:stress:d} $4 l_1/c_0 = 187.8$ $\mu$s \subref{fig:longitudinal_impact_of_two_bars:stress:e} $5 l_1/c_0 = 187.8$ $\mu$s. Two numerical simulations are considered with OBS (red line) and EBS (blue dashed line). The analytical solutions are also plotted with black dotted line. EBS yields better accuracy over OBS and converges to the analytical solutions.}
	\label{fig:longitudinal_impact_of_two_bars:stress}
\end{figure}

Both bars possess identical material properties, with Young's modulus $E=50.5$ GPa, cross section $A=1$ m\textsuperscript{2} and mass density $\rho=2783$ kg/m\textsuperscript{3}. All the points on the left bar have the same initial velocity, i.e. $\dot{u}_1=v_0=1$ m/s, while the bar on the right is stationary $\dot{u}_2=0$ m/s. The two bars are unstressed at $t=0$. The solution of this problem can be derived from the one-dimensional wave equation as show in \FIG{fig:longitudinal_impact_of_two_bars:geo_ibcs}. The simulation is performed with a grid spacing $\Delta h = 0.390625$ mm and the material points are initially located at the Gauss positions of their corresponding grid cells. The cell density is $4$ material points per grid cell while extra material points are added at the bars edges to track their boundary. A total of $6148$ material points was used. Time step is chosen to be $\Delta t = 0.00917$ $\mu$s and satisfies the Courant– Friedrichs– Lewy condition. The value of the normal penalty parameter is $\omega^{nor} = 129280$ N/mm\textsuperscript{3}. To alleviate any stress “noise” during impact of the two bars, the proposed EBS-based MPM is utilised with occupation parameter $C_c=0.75$. Simulation is also performed with OBS to highlight the merits of the proposed method. The one-dimensional Eulerian grid is formed with $x^{min} = 0$ m, $x^{max} = 1$ m. The results from the simulation are shown in \FIG{fig:longitudinal_impact_of_two_bars:vel} and \FIG{fig:longitudinal_impact_of_two_bars:stress} for the bars' velocity and stress over $5$ time instances, $t \in [l_1/c_0, 2 l_1/c_0, 3 l_1/c_0, 4 l_1/c_0, 5 l_1/c_0] = [46.95,	93.9,	140.85,	187.8,	234.75]$ $\mu$s.

The two objects undergo contact until $t = 4 l_1 / c_0$, which is the point in time at which the reflected wave in the second object reaches the point of contact. Given the stress-free nature of the first object, the wave reflects back instead of entering the first object in accordance with the stress-free boundary condition. This is illustrated in \FIG{fig:longitudinal_impact_of_two_bars:vel:c} and \FIG{fig:longitudinal_impact_of_two_bars:stress:c} where both the mean value of the velocity of the stress are zero. Following this, the objects lose contact. The final velocity of the first object is $\dot{u}_1=0$, whereas that of the second object is $\dot{u}_2=v_0/2$. This is shown in \FIG{fig:longitudinal_impact_of_two_bars:vel:d}. Notably, as depicted in \FIG{fig:longitudinal_impact_of_two_bars:vel:e}, the second object continues to undergo oscillations as a result of the travelling stress wave, while the first object remains stationary. 

Velocity and stress numerical solutions also result good agreement with the analytical solutions provided in \cite{Wriggers2006ComputationalMechanics}. The numerical solutions with OBS are also plotted in \FIG{fig:longitudinal_impact_of_two_bars:vel} and \FIG{fig:longitudinal_impact_of_two_bars:stress} with red colour. By comparing the two numerical solutions, it is pronounced that OBS fails to accurately capture the velocity and stress distribution along the bars, yielding in severe errors mainly in the location of the contact point. 

\begin{figure}[htbp]
	\centering
	\begin{tabular}{cc}
            \multicolumn{2}{c}{\includegraphics[scale=0.85]{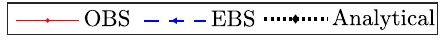}}\\
		  \subfloat[\label{fig:longitudinal_impact_of_two_bars:stress_parametric:a}]{
			\includegraphics[scale=0.85]{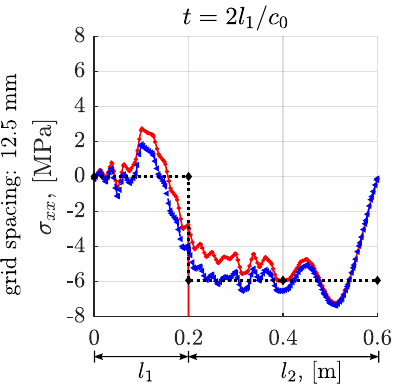}} &
            \subfloat[\label{fig:longitudinal_impact_of_two_bars:stress_parametric:b}]{
			\includegraphics[scale=0.85]{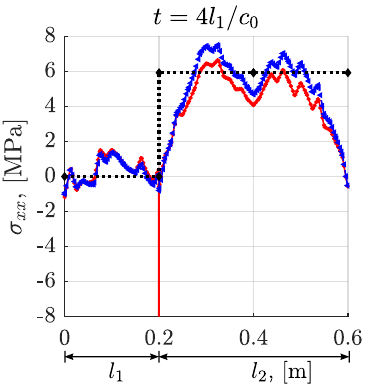}} \\ 
            \subfloat[\label{fig:longitudinal_impact_of_two_bars:stress_parametric:c}]{
			\includegraphics[scale=0.85]{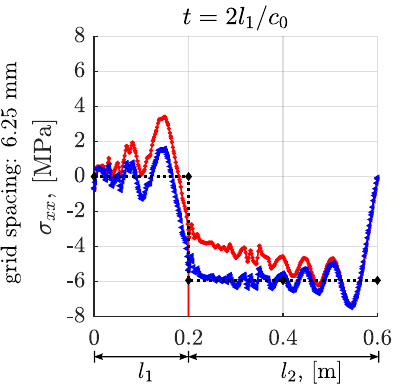}} &
            \subfloat[\label{fig:longitudinal_impact_of_two_bars:stress_parametric:d}]{
			\includegraphics[scale=0.85]{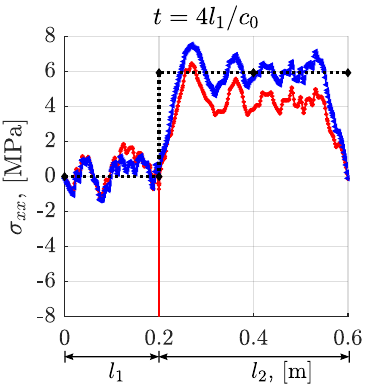}} \\
            \subfloat[ \label{fig:longitudinal_impact_of_two_bars:stress_parametric:e}]{
			\includegraphics[scale=0.85]{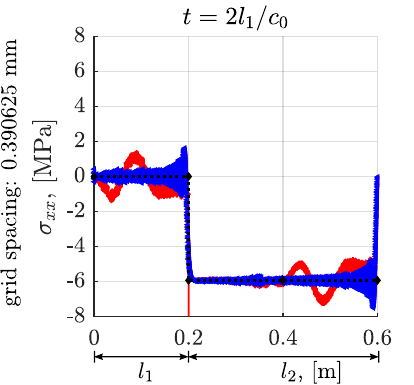}} & 
            \subfloat[ \label{fig:longitudinal_impact_of_two_bars:stress_parametric:f}]{
			\includegraphics[scale=0.85]{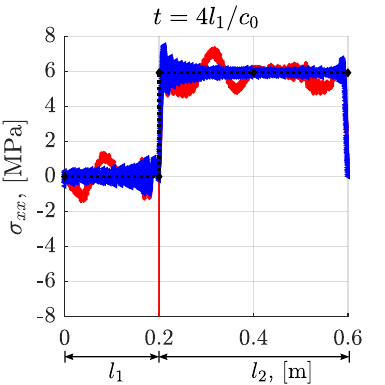}}   
	\end{tabular}
	\caption[]{Longitudinal impact of two bars. Stress over bars' length for time instances, $t \in [2 l_1/c_0, 4 l_1/c_0] = [93.9,	187.8]$ $\mu$s and three grid spacings, $\Delta h \in [12.5, 6.25, 0.390625]$ mm. Two numerical simulations are considered with OBS (red line) and EBS (blue dashed line). The analytical solutions are also plotted with black dotted line. EBS yields better accuracy over OBS and converges to the analytical solutions faster.}
	\label{fig:longitudinal_impact_of_two_bars:stress_parametric}
\end{figure}

As a further investigation within this example, a parametric study is performed for the grid spacing of the Eulerian grid as shown in \FIG{fig:longitudinal_impact_of_two_bars:stress_parametric}. Two additional coarser grid spacings are considered, i.e. $12.5$ mm and $6.25$ mm. The time step is then chosen as $0.2934$ $\mu$s and $0.1467$ $\mu$s, respectively. Similarly, the number of materials points used are $196$ and $388$, respectively. The value of the normal penalty parameter is also adjusted at $4040$ N/mm\textsuperscript{3} and $8080$ N/mm\textsuperscript{3}. \FIG{fig:longitudinal_impact_of_two_bars:stress_parametric} presents the bars' stress over their length for two time instances, i.e. $t = 2 l_1 / c_0$ and $t = 4 l_1 / c_0$. In this, it is shown that EBS yields better accuracy over OBS and converges to the analytical solutions faster.

\subsection{Hertz’s contact benchmark problems} \label{sec:hertz_contact}

These examples illustrate the proposed method's ability to evaluate contact areas using Hertz contact analyses.

\subsubsection{Scenario 1: Hertz contact between two disks} \label{subsec:hertz_contact_Guilkey}

\begin{figure}[htbp]
    \centering
    \includegraphics[scale=1.0]{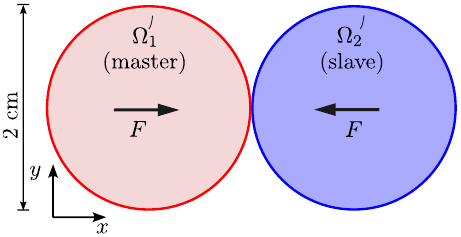}
	\caption{Hertz contact between two disks. Scenario 1. Geometry and boundary conditions.}
	\label{fig:hertz_contact_Guilkey_geo_bcs}
\end{figure}

In the first Hertz contact problem, we examine two contacting disks as depicted in \FIG{fig:hertz_contact_Guilkey_geo_bcs}. Both disks have a radius of $R = 1$ cm, a Young's modulus of $9.453 \cdot 10^{10}$ Pa, and a Poisson's ratio of $0.07$. The disks experience body forces $F$ that compel them to approach each other. This benchmark has previously been investigated by \cite{Guilkey2021AMethod} using a hybrid penalty and grid-based contact method for MPM. Their results provide a basis for comparison with our proposed approach. In particular, \cite{Guilkey2021AMethod} determined the force-contact length relationship for applied force $F$ across three grid spacings, $\Delta h$: $R/10$, $R/20$, and $R/40$. They explored two orientations of the line through the disk centres relative to the X-axis: (i) $0$ degrees and (ii) $45$ degrees.

We conduct comparable simulations to those in \cite{Guilkey2021AMethod}. The bars are stress-free at $t=0$, and a quasi-static response is achieved by incrementally raising body forces $F$ over $400$ $\mu$s. For consistency with \cite{Guilkey2021AMethod}, we use $100$ master segments for the boundary of each disk in the $0$ degrees configuration and $200$ segments for the $45$ degrees configuration. The simulated contact length is determined by aggregating penetrated master segments.

\begin{figure}[htbp]
	\centering
	\begin{tabular}{c}
            \multicolumn{1}{c}{\includegraphics[scale=1.0]{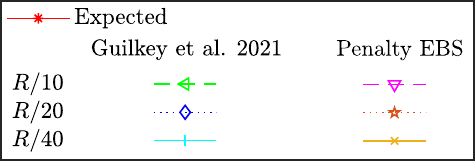}}\\
		  \subfloat[\label{fig:hertz_contact_Guilkey:force_contact:a}]{
			\includegraphics[scale=1.0]{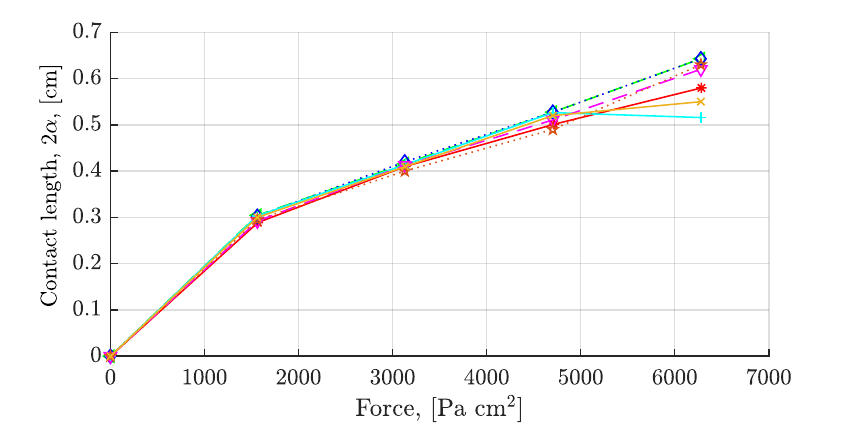}} \\
            \subfloat[\label{fig:hertz_contact_Guilkey:force_contact:b}]{
			\includegraphics[scale=1.0]{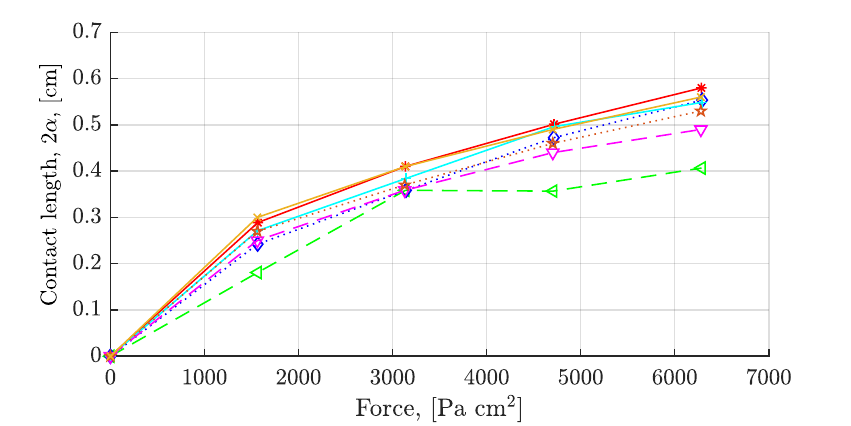}}   
	\end{tabular}
	\caption[]{Hertz contact between two disks. Scenario 1. Contact length analysis for two disks pushed together by body force. Numerical results at three grid resolutions are compared to Hertz contact analysis results and the simulations in \cite{Guilkey2021AMethod}. \subref{fig:hertz_contact_Guilkey:force_contact:a} $0$ degrees and $100$ master segments, and \subref{fig:hertz_contact_Guilkey:force_contact:b} $45$ degrees and $200$ master segments.}
    \label{fig:hertz_contact_Guilkey:force_contact}
\end{figure}

Simulation results for the $0$ and $45$ degrees configurations are displayed in \FIG{fig:hertz_contact_Guilkey:force_contact:a} and \FIG{fig:hertz_contact_Guilkey:force_contact:b}, respectively. Analytical expressions of the force-contact length relationship, detailed in \cite{Guilkey2021AMethod, Johnson1985ContactMechanics}, are also included in the figures for comparison with the numerical simulations. The analytical relationship between the contact force $F$ and the contact region (semi-contact-width) $a$ is
\begin{equation}
\alpha = \sqrt{\frac{4 F R^{*}}{\pi E^{*} L}}
\label{eqn:hertz_contact:scenario_1_contact_surface}
\end{equation}
where $\frac{1}{R^{*}} = \frac{1}{R_1} + \frac{1}{R_2}$ and $\frac{1}{E^{*}} = \frac{1-\nu^{2}_1}{E_1} + \frac{1-\nu^{2}_2}{E_2}$ are the effective radius and elastic modulus, respectively and $L$ the length of the cylinder assumed to be $1$ m.

All simulations using the proposed method align well with expected and computed results from \cite{Guilkey2021AMethod}. As noted in \cite{Guilkey2021AMethod} and depicted in \FIG{fig:hertz_contact_Guilkey:force_contact:a}, deviations from expected results are more noticeable in the final loading step. The responses of cases $R/10$ and $R/20$ are similar, while the $R/40$ prediction underestimates the contact area, albeit providing more accurate estimates than those reported in \cite{Guilkey2021AMethod}.

\begin{figure}[htbp]
	\centering
	\begin{tabular}{cc}
            \multicolumn{2}{c}{\includegraphics[scale=1.0]{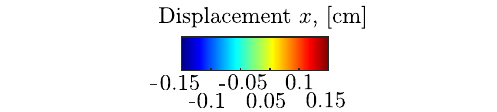}}\\
		  \subfloat[\label{fig:hertz_contact_Guilkey:force_contact_visual:a}]{
			\includegraphics[scale=1.00]{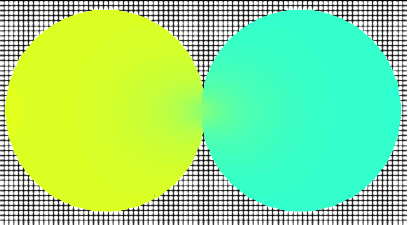}} &
            \subfloat[\label{fig:hertz_contact_Guilkey:force_contact_visual:b}]{
			\includegraphics[scale=1.00]{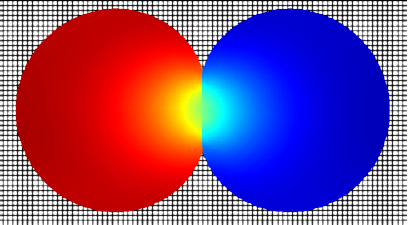}}   
	\end{tabular}
	\caption[]{Hertz contact between two disks. Scenario 1. Contact length analysis for two disks pushed together by body force. Penalty EBS results. Displacement $x$ for grid spacing $R/20$, $0$ degrees and $100$ master segments for the domain boundary. \subref{fig:hertz_contact_Guilkey:force_contact_visual:a} first loading step and  \subref{fig:hertz_contact_Guilkey:force_contact_visual:b} last loading step. These loading steps correspond to the markers shown in \FIG{fig:hertz_contact_Guilkey:force_contact:a}.}
    \label{fig:hertz_contact_Guilkey:force_contact_visual}
\end{figure}

In \FIG{fig:hertz_contact_Guilkey:force_contact:b}, excellent agreement with the expected result is also observed for the $45$ degrees case. Notably, the simulation estimates show better alignment with the expected result when the grid spacing resolution is finer, specifically with $R/40$. The increased resolution of the domain boundary, achieved with $200$ master segments, further improves prediction accuracy, as also highlighted in \cite{Guilkey2021AMethod}. Importantly, the results obtained with the current method exhibit better agreement with the expected outcome compared to those reported in \cite{Guilkey2021AMethod}. A visual representation of the results obtained with the proposed penalty EBS are shown in \FIG{fig:hertz_contact_Guilkey:force_contact_visual} for the first and last loading steps.

\subsubsection{Scenario 2: Hertz contact between a demi-sphere and a rigid plane} \label{subsec:hertz_contact_Sun}

\begin{figure}[htbp]
    \centering
    \includegraphics[scale=1.0]{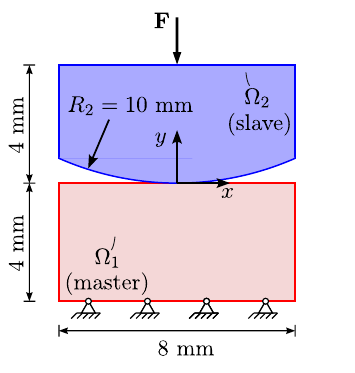}
	\caption{Hertz contact between a demi-sphere and a rigid plane. Scenario $2$. Geometry and boundary conditions.}
	\label{fig:hertz_contact_geo_bcs}
\end{figure}

Scenario $2$ explores Hertz's contact benchmark, involving a demi-sphere and a rigid plane. This aims to validate the proposed MPM implementation for non-flat contact surfaces by comparing the analytical solution of the pressure distribution and the contact area size. Previous studies \cite{Xiao2021DP-MPM:Fragmentation} are referenced for discussion. The problem's configuration and boundary conditions are illustrated in \FIG{fig:hertz_contact_geo_bcs}.

For the slave body, Young’s modulus is $E_{2} = 10$ GPa, and Poisson’s ratio is $\nu_{2} = 0$. The master body is rigid. Three Eulerian grid spacings are considered: $\Delta h = 0.2$ mm, $\Delta h = 0.1$ mm, and $\Delta h = 0.05$ mm. The simulation duration is $175$ $\mu$s. An occupation parameter of $0.50$ is chosen for EBS generation, and the friction coefficient is $\mu_f = 0$. The two bodies are unstressed at $t=0$. An external load of $F = 156.7$ N is gradually applied throughout the simulation to mimic conditions approaching a steady state.

The analytical solution of the contact pressure distribution along the contact surface is evaluated from the relation below
\begin{equation}
\sigma_{yy} = - \sigma^{max}_{yy} \sqrt{1 - \left( \frac{s}{b} \right)^2}, \text{   for   } 0 \leq s \leq b
\label{eqn:hertz_contact:contact_pressure}
\end{equation}
where $\sigma^{max}_{yy} = \frac{2F}{\pi b}$ and $b=2\sqrt{\frac{\left( 2 F R \right)}{\pi E^{\prime}}}$. The effective Young's modulus is defined as $E^{\prime} = \frac{2}{E^{\prime}}=\frac{1-\nu_1^2}{E_1}+\frac{1-\nu_2^2}{E_2}$ while $R=\frac{R_1 R_2}{R_1 + R_2}$ is the equivalent body radius, where $R_1$, $R_2$ are the radii of the two contact surfaces, respectively. The radius of the master body is infinite.

\begin{figure}[htbp]
	\centering
	\begin{tabular}{cc}
            \multicolumn{2}{c}{\includegraphics[scale=1.0]{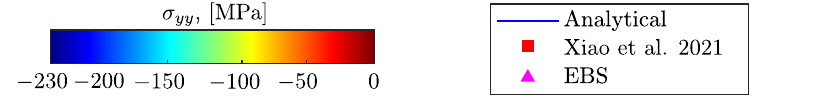}}\\
		  \subfloat[\label{fig:hertz_contact_xiao_ebs_syy:a}]{
			\includegraphics[scale=1.0]{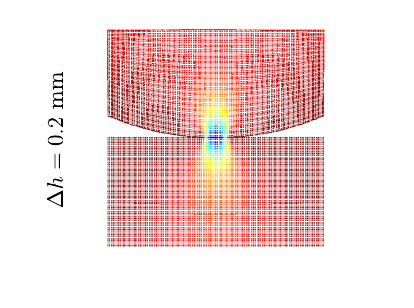}} &
            \subfloat[\label{fig:hertz_contact_xiao_ebs_syy:b}]{
			\includegraphics[scale=1.0]{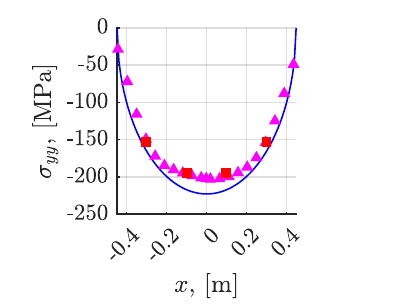}} \\       \subfloat[\label{fig:hertz_contact_xiao_ebs_syy:c}]{
			\includegraphics[scale=1.0]{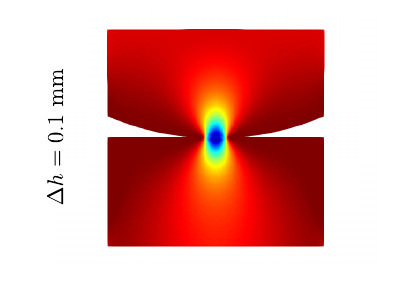}} &
            \subfloat[\label{fig:hertz_contact_xiao_ebs_syy:d}]{
			\includegraphics[scale=1.0]{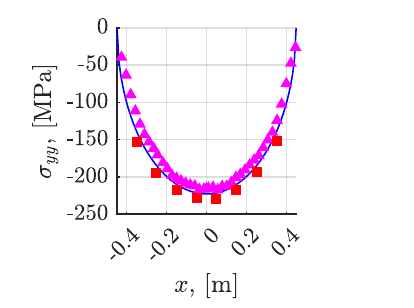}} \\
            \subfloat[\label{fig:hertz_contact_xiao_ebs_syy:e}]{
			\includegraphics[scale=1.0]{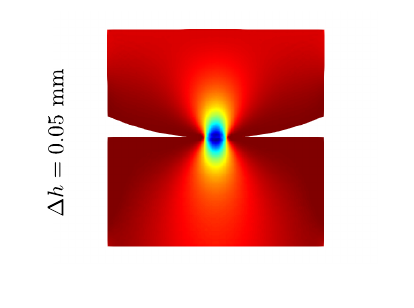}} &
            \subfloat[\label{fig:hertz_contact_xiao_ebs_syy:f}]{
			\includegraphics[scale=1.0]{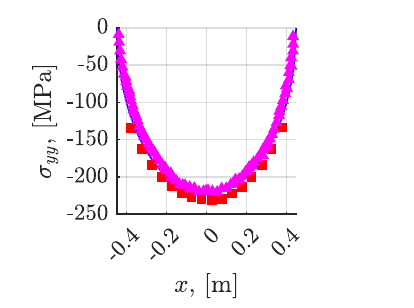}}   
	\end{tabular}
	\caption[]{Hertz contact between a demi-sphere and a rigid plane. Scenario $2$. Results of $\sigma_{yy}$ with the proposed EBS-based MPM and three grid spacings, i.e. \subref{fig:hertz_contact_xiao_ebs_syy:a}, \subref{fig:hertz_contact_xiao_ebs_syy:b} $\Delta h = 0.2$ mm, \subref{fig:hertz_contact_xiao_ebs_syy:c}, \subref{fig:hertz_contact_xiao_ebs_syy:d} $\Delta h = 0.1$ mm and \subref{fig:hertz_contact_xiao_ebs_syy:e}, \subref{fig:hertz_contact_xiao_ebs_syy:f} $\Delta h = 0.05$ mm. The simulation results imply that, as the mesh undergoes iterative refinement, the contact pressure progressively nears and aligns with the analytical solution. This is also reported in the findings of \cite{Xiao2021DP-MPM:Fragmentation}, which are marked in the graphs.}
    \label{fig:hertz_contact_xiao_ebs_syy}
\end{figure}

The $\sigma_{yy}$ results from the EBS-based MPM are depicted in \FIG{fig:hertz_contact_xiao_ebs_syy} for the three grid spacings. The results suggest a convergence of contact pressure to the analytical solution with increasing mesh refinement, consistent with \cite{Xiao2021DP-MPM:Fragmentation}.

\begin{figure}[htbp]
    \centering
    \includegraphics[scale=1.0]{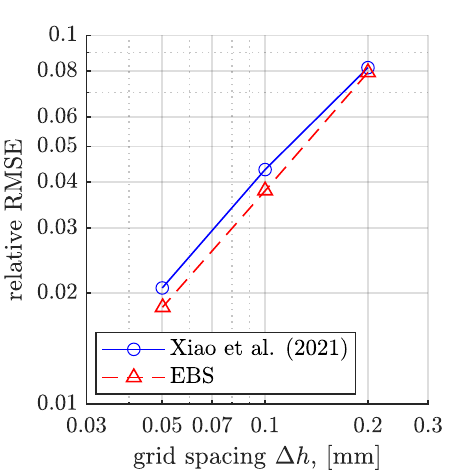}
	\caption{Hertz contact between a demi-sphere and a rigid plane. Scenario $2$. Convergence of the relative root mean square error (RMSE) (see \EQ{eqn:hertz_contact:error}) for simulated contact pressures using various grid spacings and EBS-based MPM. Computed error is comparable with the results reported in \cite{Xiao2021DP-MPM:Fragmentation}.}	\label{fig:hertz_contact_xiao_ebs_error}
\end{figure}

To further assess the proposed approach's accuracy, \FIG{fig:hertz_contact_xiao_ebs_error} illustrates the normalized relative Root Mean Square Error (RMSE) for pressure values across contact surface material points. The relative RMSE is formulated as
\begin{equation}
\text{relative RMSE} = \frac{\sqrt{ \frac{\sum\limits_{p = 1}^{{\mathscr{N}_{p}}} { \left( {\sigma}^{\text{numerical}}_{yy} - {\sigma}^{\text{analytical}}_{yy} \right)^2}}{\mathscr{N}_{p}}}}{\sigma^{max}_{yy}}
\label{eqn:hertz_contact:error}
\end{equation}
where $\mathscr{N}_p$ is the total number of material points at the contact surface $\partial \Omega_2$. As anticipated, the relative Root Mean Square Error (RMSE) diminishes as the refinement of the Eulerian grid takes place.

\begin{figure}[htbp]
	\centering
	\begin{tabular}{ccc}
            \multicolumn{3}{c}{\includegraphics[scale=1.0]{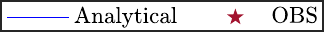}}\\
		  \subfloat[\label{fig:hertz_contact_obs_syy:a}]{
			\includegraphics[scale=1.0]{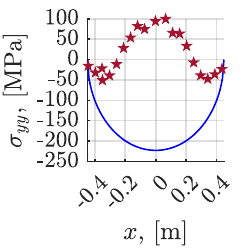}} &
            \subfloat[\label{fig:hertz_contact_obs_syy:b}]{
			\includegraphics[scale=1.0]{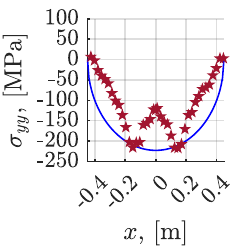}} &
            \subfloat[\label{fig:hertz_contact_obs_syy:c}]{
			\includegraphics[scale=1.0]{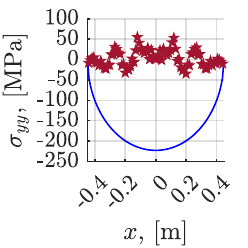}}
	\end{tabular}
	\caption[]{Hertz contact between a demi-sphere and a rigid plane. Scenario $2$. Results of $\sigma_{yy}$ with the OBS-based MPM and three grid spacings, i.e. \subref{fig:hertz_contact_obs_syy:a} $\Delta h = 0.2$ mm, \subref{fig:hertz_contact_obs_syy:b} $\Delta h = 0.1$ mm and \subref{fig:hertz_contact_obs_syy:c} $\Delta h = 0.05$ mm. As expected, OBS exhibits a lack of convergence with the analytical solution even when the Eulerian grid is fine.}
    \label{fig:hertz_contact_obs_syy}
\end{figure}

For a comprehensive view, \FIG{fig:hertz_contact_obs_syy} presents the corresponding OBS results across the three grid spacings. Similar to the one-dimensional problems discussed earlier, it's clear that OBS struggles to attain accurate solutions for contact stresses, even with a fine Eulerian grid.

\subsection{Impact of two compressible Neo-Hookean cylinders} \label{sec:imp_elast_rings}

\begin{figure}[htbp]
    \centering
    \includegraphics[scale=1.0]{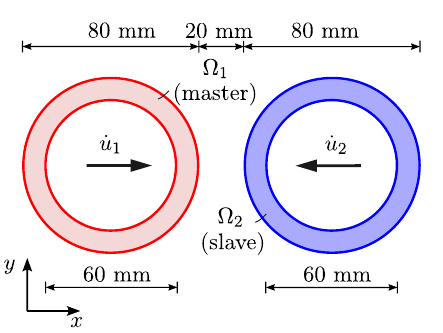}
	\caption{Impact of two compressible Neo-Hookean cylinders. Geometry and boundary conditions.}
	\label{fig:imp_elast_rings_geo_bcs}
\end{figure}

This example examines the collision of two hollow elastic cylinders as shown in \FIG{fig:imp_elast_rings_geo_bcs}. This benchmark has  been extensively investigated in the litarature, see, e.g.,  \cite{Huang2011ContactSimulation} where discrete velocity fields \cite{Bardenhagen2001AnMaterial} have been employed and most recently by \cite{deVaucorbeil2021ModellingMethod} using a TLMPM formulation. The aim of this example is demonstrate that the proposed method results in no early contact errors. Both cylinders are considered purely elastic, and hence energy conservation is also examined.

The material is a compressible Neo-Hookean which the stored energy density is defined as \cite{Bonet2008NonlinearAnalysis}
\begin{equation}
\psi = \frac{\mu}{2} \left( I_c - d \right) - \mu \ln J + \frac{\lambda}{2} (\ln J)^2
\label{eqn:imp_elast_rings_stored_energy}
\end{equation}
where $I_c = \Tr{(\mathbf{C})}$ is the first principal invariant of the right Cauchy-Green strain tensor $\mathbf{C} = \mathbf{F}^{T} \mathbf{F}$ and $J = \det{(\mathbf{F})}$. The symbol $\Tr$ denotes the trace of the tensor $\mathbf{C}$ while $d =2$ corresponds to the dimension of the problem. The material properties are expressed with the Lam\'{e} constants, i.e. $\mu = 26.1$ MPa and $\lambda = 104.4$ MPa. The Cauchy stress tensor derives from \EQ{eqn:imp_elast_rings_stored_energy} and is expressed as 
\begin{equation}
\sigma = \frac{1}{J} \left( \mu (\mathbf{F} \mathbf{F}^{T} - \mathbf{I}) + \lambda \ln(J \mathbf{I}) \right)
\label{eqn:imp_elast_rings_stress}
\end{equation}
where $\mathbf{I}$ is the second-order identity tensor.

The mass density is $\rho = 1010$ kg/m\textsuperscript{3} and the grid cell size is considered to be $\Delta h = 1.25$ mm. The magnitude of the cylinders' initial velocity is $0.03$ mm/$\mu$s. Frictional contact is not employed in this configuration, i.e. $\mu_f = 0$. Both solid bodies are unstressed at $t=0$.

The kinetic, stored and total energy for each time step for all material points are evaluated as
\begin{equation}
\mathcal{K} = \sum\limits_{p = 1}^{{n_p}} {\left(\frac{1}{2} (\dot{u}_{p_x}^2 + \dot{u}_{p_y}^2) M_p \right)} \text{  ,  } \mathcal{W} = \sum\limits_{p = 1}^{{n_p}} {\left( \psi_p \Omega_p \right)} \text{  ,  } \mathcal{T} = \mathcal{K} + \mathcal{W},
\label{eqn:imp_elast_rings_kin_eng}
\end{equation}
respectively.

\begin{figure}[htbp]
	\centering
	\begin{tabular}{cc}
            \multicolumn{2}{c}{\includegraphics[scale=1.0]{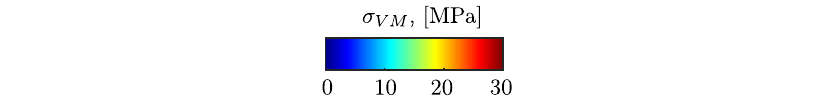}}\\
		  \subfloat[\label{fig:imp_elast_rings_ebs_svn:a}]{
			\includegraphics[scale=1.0]{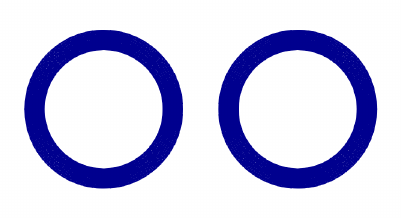}} &
            \subfloat[\label{fig:imp_elast_rings_ebs_svn:b}]{
			\includegraphics[scale=1.0]{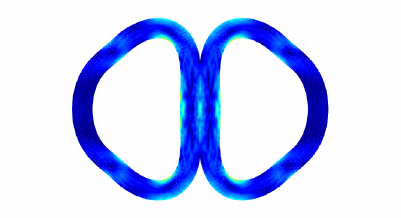}} \\       \subfloat[\label{fig:imp_elast_rings_ebs_svn:c}]{
			\includegraphics[scale=1.0]{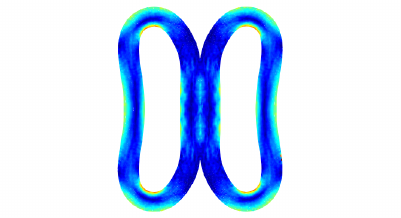}} &
            \subfloat[\label{fig:imp_elast_rings_ebs_svn:d}]{
			\includegraphics[scale=1.0]{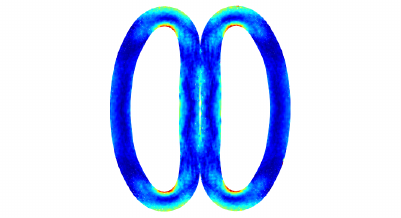}} \\
            \subfloat[\label{fig:imp_elast_rings_ebs_svn:e}]{
			\includegraphics[scale=1.0]{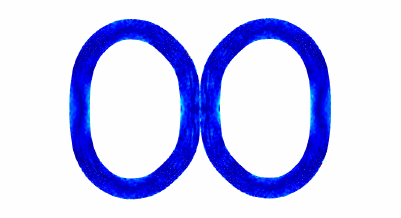}} &
            \subfloat[\label{fig:imp_elast_rings_ebs_svn:f}]{
			\includegraphics[scale=1.0]{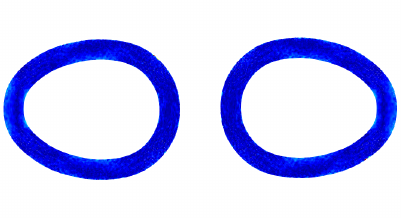}}   
	\end{tabular}
	\caption[]{Impact of two compressible Neo-Hookean cylinders. Deformation evolution for time steps \subref{fig:imp_elast_rings_ebs_svn:a} $t=0$ ms, \subref{fig:imp_elast_rings_ebs_svn:b} $t=0.87$ ms, \subref{fig:imp_elast_rings_ebs_svn:c} $t=1.59$ ms, \subref{fig:imp_elast_rings_ebs_svn:d} $t= 2.33$ ms, \subref{fig:imp_elast_rings_ebs_svn:e} $t=3.09$ ms and \subref{fig:imp_elast_rings_ebs_svn:f} $t=3.97$ ms. The colours correspond to the magnitude of the Von Mises stress. Results obtained with penalty EBS with $C_c = 0.50$. It is evident from the above figures that the proposed method results in no early contact error when the two cylinders make contact.}
    \label{fig:imp_elast_rings_ebs_svn}
\end{figure}

\begin{figure}[htbp]
	\centering
	\begin{tabular}{ccc}
            \multicolumn{3}{c}{\includegraphics[scale=1.0]{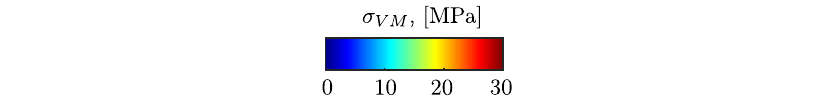}}\\
		  \subfloat[\label{fig:imp_elast_rings_ebs_TLMPM:a}]{
			\includegraphics[scale=1.0]{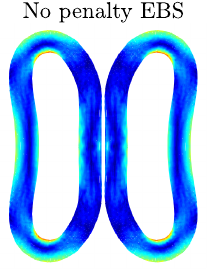}} &
            \subfloat[\label{fig:imp_elast_rings_ebs_TLMPM:b}]{
			\includegraphics[scale=1.0]{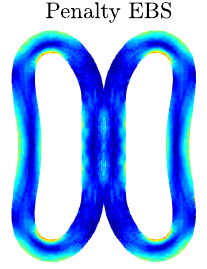}} &
            \subfloat[\label{fig:imp_elast_rings_ebs_TLMPM:c}]{
			\includegraphics[scale=1.0]{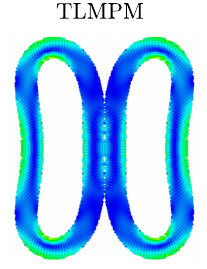}}
	\end{tabular}
	\caption[]{Impact of two compressible Neo-Hookean cylinders. Von Mises stress at time $t=1.59$ ms for \subref{fig:imp_elast_rings_ebs_TLMPM:a} no penalty EBS, \subref{fig:imp_elast_rings_ebs_TLMPM:b} proposed penalty EBS and \subref{fig:imp_elast_rings_ebs_TLMPM:c} TLMPM \cite{deVaucorbeil2021ModellingMethod}. The no penalty EBS results in early contact error whereas both the proposed penalty EBS and TLMPM have successfully eliminated such errors. However, the penalty EBS demonstrate better energy conservation compared to TLMPM as shown in \FIG{fig:imp_elast_rings_engs}. The TLMPM results shown here derived by utilising hat weighting functions which as per \cite{deVaucorbeil2021ModellingMethod} yield better energy conservation.}
    \label{fig:imp_elast_rings_ebs_TLMPM}
\end{figure}

\begin{figure}[htbp]
	\centering
	\begin{tabular}{cc}   
   \multicolumn{2}{c}{\includegraphics[scale=1.0]{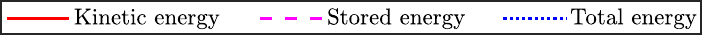}}\\
		\subfloat[\label{fig:imp_elast_rings_engs:ebs}]{
			\includegraphics[scale=1.0]{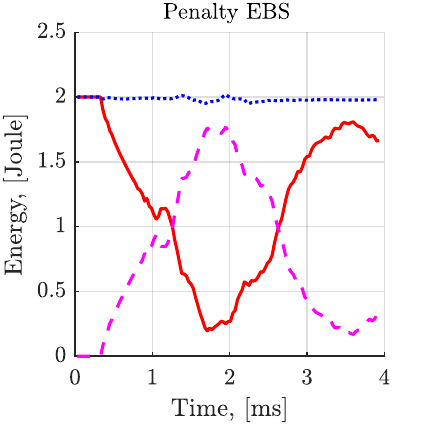}} &
            \subfloat[\label{fig:imp_elast_rings_engs:tlmpm}]{
			\includegraphics[scale=1.0]{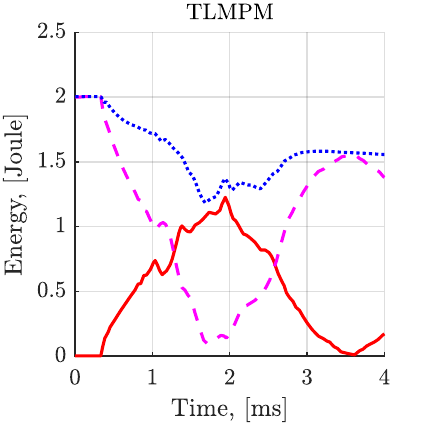}}
	\end{tabular}
	\caption[]{Impact of two compressible Neo-Hookean cylinders. Kinetic, stored and total energies for \subref{fig:imp_elast_rings_engs:ebs} penalty EBS and \subref{fig:imp_elast_rings_engs:tlmpm} TLMPM as proposed in \cite{deVaucorbeil2021ModellingMethod}. The TLMPM results presented here are obtained using hat weighting functions, which, according to \cite{deVaucorbeil2021ModellingMethod}, result in enhanced energy conservation. The proposed penalty EBS method yields better energy conservation results.}
	\label{fig:imp_elast_rings_engs}
\end{figure}

The evolution of deformation over time is depicted across various time steps in \FIG{fig:imp_elast_rings_ebs_svn}, employing the proposed penalty EBS method. \FIG{fig:imp_elast_rings_ebs_TLMPM} also contrasts the proposed method with the no penalty EBS and TLMPM formulation as outlined in \cite{deVaucorbeil2021ModellingMethod}. The no penalty EBS results were obtained by considering a common discrete field for both solid bodies. In \FIG{fig:imp_elast_rings_ebs_TLMPM}, it becomes apparent that when a penalty method is not used the MPM implementation leads to early contact errors as the two solid bodies make contact. In comparison, both the proposed penalty EBS and the TLMPM variants have successfully eliminated such errors. However, as illustrated in \FIG{fig:imp_elast_rings_engs}, the results obtained with the proposed method demonstrate better energy conservation compared to TLMPM. The TLMPM results depicted in \FIG{fig:imp_elast_rings_ebs_TLMPM:c} and \FIG{fig:imp_elast_rings_engs:tlmpm} are obtained using hat weighting functions, which, as suggested in \cite{deVaucorbeil2021ModellingMethod}, contribute to enhanced energy conservation.

\subsection{Stress wave in a granular material} \label{sec:stress_wave_gran_mat}

The previous examples focused solely on contact between two solid entities, i.e. two separate discrete fields, whereas the proposed algorithm can handle multiple contacts. In this instance, we explore two significant stress wave propagation scenarios in granular media, as discussed in \cite{Bardenhagen2001AnMaterial}.

\subsubsection{Scenario 1} \label{sec:stress_wave_gran_mat:scenario_1}

In the first scenario, four identical collinear disks are impacted by a right-travelling striker with an initial velocity of $0.0056$ mm/$\mu$s. The geometry and boundary conditions are summarised visually in \FIG{eqn:stress_wave_gran_mat:stress_fringes}. The problem is divided into $5$ separate discrete fields, forming a total of $10$ discrete pairs, as illustrated in the figure. Both the disks and the impactor are made of linear elastic material. The material properties for the impactor are $E = 17485.71$ MPa, $\nu = 0.214$, $\rho = 190000$ kg/m\textsuperscript{3}, and for each disk, they are $E = 174857.14$ MPa, $\nu = 0.214$, $\rho = 1900$ kg/m\textsuperscript{3}. Frictional contact is not employed in this configuration ($\mu_f = 0$) as there is no sliding among the disks.

The grid cell density is set at $16$ material points per grid cell element for both the impactor and the disks. Three grid cell size cases are considered: (i) $\Delta h = 5$ mm, (ii) $0.025$ $\mu$s, and (iii) $0.0125$ $\mu$s.

Qualitative comparisons between EBS results and experimental findings from \cite{Bardenhagen2001AnMaterial} are conducted. The stress distribution in the disks, assessed through photo-elasticity in the experiments, is represented by dark fringes at contours of constant maximum difference in principal stresses. The simulations produce fringes using the formula:
\begin{equation}
\begin{aligned}
\sigma_f = 1 - \sin^2 \left( k_f (\sigma_1 - \sigma_2) \right)
\end{aligned}
\label{eqn:stress_wave_gran_mat:stress_fringes}
\end{equation}
where $k_f = \pi / 0.07$ GPa. The disparity in in-plane principal stress is expressed as:
\begin{equation}
\begin{aligned}
\sigma_1 - \sigma_2 = \sqrt{4 \tau^{2}_{xy} + (\sigma_{xx} - \sigma_{yy})^2 }
\end{aligned}
\label{eqn:stress_wave_gran_mat:princ_stress_diff}
\end{equation}

\begin{figure}[htbp]
    \centering
    \includegraphics[scale=0.85]{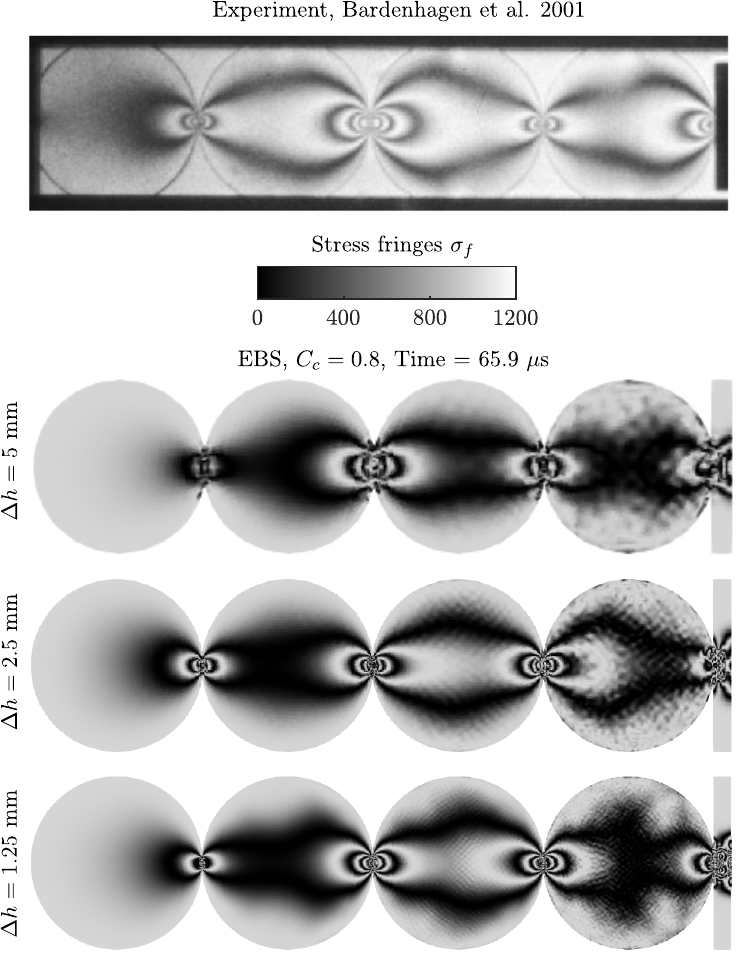}
	\caption{Stress wave in granular media: Scenario 1. The stress fringe pattern derived through EBS at time instance $65.9$ $\mu$s and the three grid spacing cases. The findings obtained from EBS yield good agreement with the experimental data from \cite{Bardenhagen2001AnMaterial}.}    \label{fig:stress_wave_gran_mat_1:experiment}
\end{figure}

The stress fringe pattern derived through EBS is presented in \FIG{fig:stress_wave_gran_mat_1:experiment} for a time instance of $80$ $\mu$s, using an occupation parameter at $C_c = 0.8$. Results from all grid spacings demonstrate a similar fringe pattern, and convergence is achieved for fine grid spacings at $2.5$ mm and $1.25$ mm. In \FIG{fig:stress_wave_gran_mat_1:experiment}, the EBS-derived fringe pattern is compared against the experimental observations from \cite{Bardenhagen2001AnMaterial}, showing good agreement.

\begin{figure}[htbp]
    \centering
    \includegraphics[scale=1.0]{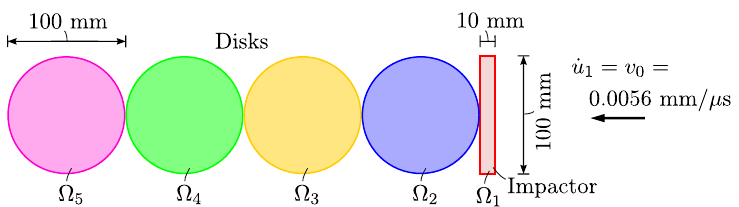}
	\caption{Stress wave in granular media: Scenario 1. Geometry and boundary conditions.}
    \label{fig:stress_wave_gran_mat_1:geo_ibcs}
\end{figure}

\subsubsection{Scenario 2} \label{sec:stress_wave_gran_mat:scenario_2}

\begin{figure}[htbp]
    \centering
    \includegraphics[scale=1.0]{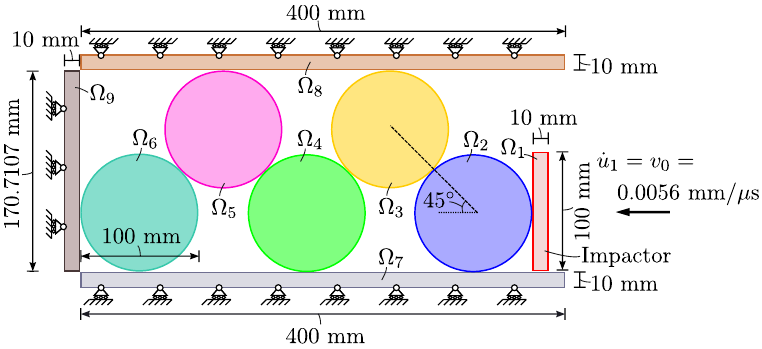}
	\caption{Stress wave in granular media: Scenario 2. Geometry and boundary conditions.}
    \label{fig:stress_wave_gran_mat_2:geo_ibcs}
\end{figure}

The proposed MPM contact algorithm facilitates the inclusion of frictional contacts, thereby enabling the simulation of stress wave propagation within granular media with frictional interactions amongst constituent elements. In the second scenario, five disks mutually contacting at a $45$-degree angle are subjected to impact from the right by an impactor at the established initial velocity of $0.0056$ mm/$\mu$s. To prevent disk displacement, an enclosure confines them, as depicted in \FIG{fig:stress_wave_gran_mat_2:geo_ibcs}. \FIG{fig:stress_wave_gran_mat_2:geo_ibcs} also shows that the problem is divided into $9$ discrete fields, leading to $36$ discrete pair combinations amongst them.

The disk, impactor, and their material attributes remain consistent with the four collinear disks case (scenario 1 in Section \ref{sec:stress_wave_gran_mat:scenario_1}). The material properties of the enclosure are chosen as in the disks for brevity. The problem's discretization follows the previous scenario, employing three cases for the background grid cell size. Frictional contacts govern interactions between disks and between the rightmost disk and impactor, characterised by a friction coefficient $\mu_f = 0.5$. Contacts between disks and the box remain frictionless (i.e. $\mu_f = 0$).

\begin{figure}[htbp]
    \centering
    \includegraphics[scale=0.93]{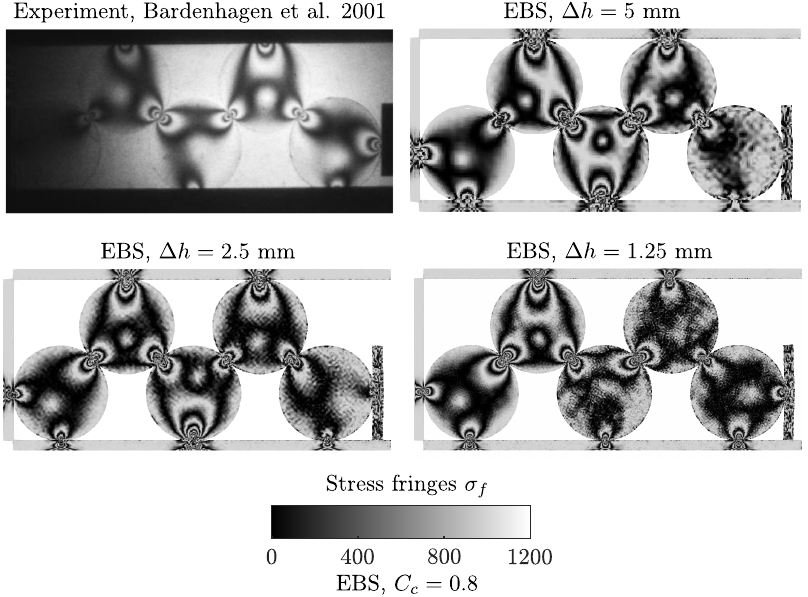}
	\caption{Stress wave in granular media: Scenario 2. The stress fringe pattern derived through EBS. Results are obtained for grid spacing at $5$ mm, $2.5$ mm and $1.25$ mm. Results from all grid spacings demonstrate similar fringe pattern and convergence is achieved for the fine grid spacing at $1.25$ mm.}   
    \label{fig:stress_wave_gran_mat_2:Dh_0.005_0.0025_0.00125m}
\end{figure}

As evident from \FIG{fig:stress_wave_gran_mat_2:Dh_0.005_0.0025_0.00125m}, a good agreement is observed between the EBS simulation results and experimental data from \cite{Bardenhagen2001AnMaterial}. The EBS results are obtained for $\Delta h = 1.25$ mm and an occupation parameter at $C_c = 0.8$. Moreover, \FIG{fig:stress_wave_gran_mat_2:Dh_0.005_0.0025_0.00125m} presents the stress fringe pattern sensitivity over grid spacing. All grid spacing cases lead to good agreement with the experiment data and convergence is achieved for the fine grid spacing at $1.25$ mm. 

\section{Conclusions} \label{sec:conclusions}

This paper introduces a novel frictional contact algorithm for material point method. It employs a penalty-based material point-to-segment contact force calculation approach. Local calculation of the normal contact force on both surfaces is achieved using a normal gap function, while the tangential contact force is computed with a tangential slip function and the Coulomb friction model. This approach improves the representation of the contact surface, preventing premature and unrealistic body contact, a limitation observed in standard MPM contact algorithms. The proposed MPM-based contact algorithm is enhanced with Extended B-Splines, addressing significant challenges such as cell-crossing and incomplete numerical integration at certain grid cells. This enrichment improves the accuracy of contact stress estimates.

The proposed method's advantages are demostrated through various numerical examples, leading to the following conclusions.
\begin{itemize}
    \item The method can be easily applied to solve frictional contact problems without significantly increasing computational complexity. Its implementation only necessitates minor modifications to the standard MPM algorithm, specifically, in the areas of basis function evaluation and handling of contact features. Furthermore, boundary tracking incurs minimal computational costs, rendering it suitable for large-scale simulations, as demonstrated in this work.
    \item The method is verified against analytical solutions and validated against experimental observations, where it yields good agreement with them.
    \item It is also compared to the standard B-Spline implementation, referred to as OBS in this work, which results in a significant improvement in the stress estimates at the contact surface.
    \item The proposed approach converges to exact solutions even with coarse meshes, unlike OBS, which necessitates a finer resolution of the Eulerian grid to reach the reference solution. To this extend, it also provides better estimates for coarser resolutions of the domain boundary against other state-of-the-art MPM contact variants.
    \item It demonstrates better energy conservation than other state-of-the-art MPM contact variants.
\end{itemize}

Future work will expand the proposed contact algorithm to address wear phenomena and self-contact issues. Ongoing developments include the creation of quasi-static and implicit time integration variants to improve the overall efficiency of the solution procedure.


\section*{Acknowledgements}
The first author is grateful to the University of Warwick for access to its high performance computing facility. The first author also acknowledges the support of the Engineering and Physical Sciences Research Council (EPSRC) - funded HetSys “Modelling of Heterogeneous Systems” consortium (EP/S022848/1). The second author acknowledges the support of the Marie Sklodowska - Curie Individual Fellowship grant “AI2AM: Artificial Intelligence driven topology optimisation of Additively Manufactured Composite Components”, No. 101021629.


\bibliography{references}

\begin{thebibliography}{10}
\expandafter\ifx\csname url\endcsname\relax
  \def\url#1{\texttt{#1}}\fi
\expandafter\ifx\csname urlprefix\endcsname\relax\def\urlprefix{URL }\fi
\expandafter\ifx\csname href\endcsname\relax
  \def\href#1#2{#2} \def\path#1{#1}\fi

\bibitem{Yang2020TheValidation}
P.~Yang, M.~Zang, H.~Zeng, X.~Guo, {The interactions between an off-road tire
  and granular terrain: GPU-based DEM-FEM simulation and experimental
  validation}, International Journal of Mechanical Sciences 179 (2020) 105634.
\newblock \href {https://doi.org/10.1016/J.IJMECSCI.2020.105634}
  {\path{doi:10.1016/J.IJMECSCI.2020.105634}}.

\bibitem{Zhang2019AnAnalysis}
W.~Zhang, E.~E. Seylabi, E.~Taciroglu, {An ABAQUS toolbox for soil-structure
  interaction analysis}, Computers and Geotechnics 114 (2019) 103143.
\newblock \href {https://doi.org/10.1016/J.COMPGEO.2019.103143}
  {\path{doi:10.1016/J.COMPGEO.2019.103143}}.

\bibitem{Chatzis2012ModelingProblem}
M.~N. Chatzis, A.~W. Smyth, {Modeling of the 3D rocking problem}, International
  Journal of Non-Linear Mechanics 47~(4) (2012) 85--98.
\newblock \href {https://doi.org/10.1016/J.IJNONLINMEC.2012.02.004}
  {\path{doi:10.1016/J.IJNONLINMEC.2012.02.004}}.

\bibitem{Marques1990ThreedimensionalForming}
M.~J. Marques, P.~A. Martins, {Three‐dimensional finite element contact
  algorithm for metal forming}, International Journal for Numerical Methods in
  Engineering 30~(7) (1990) 1341--1354.
\newblock \href {https://doi.org/10.1002/NME.1620300708}
  {\path{doi:10.1002/NME.1620300708}}.

\bibitem{Sauer2013LocalScheme}
R.~A. Sauer, {Local finite element enrichment strategies for 2D contact
  computations and a corresponding post-processing scheme}, Computational
  Mechanics 52~(2) (2013) 301--319.
\newblock \href {https://doi.org/10.1007/S00466-012-0813-8}
  {\path{doi:10.1007/S00466-012-0813-8}}.

\bibitem{Duong2019AMethod}
T.~X. Duong, L.~De~Lorenzis, R.~A. Sauer, {A segmentation-free isogeometric
  extended mortar contact method}, Computational Mechanics 63 (2019) 383--407.
\newblock \href {https://doi.org/10.1007/s00466-018-1599-0}
  {\path{doi:10.1007/s00466-018-1599-0}}.

\bibitem{Laursen2012MortarFormulations}
T.~A. Laursen, M.~A. Puso, J.~Sanders, {Mortar contact formulations for
  deformable–deformable contact: Past contributions and new extensions for
  enriched and embedded interface formulations}, Computer Methods in Applied
  Mechanics and Engineering 205-208~(1) (2012) 3--15.
\newblock \href {https://doi.org/10.1016/J.CMA.2010.09.006}
  {\path{doi:10.1016/J.CMA.2010.09.006}}.

\bibitem{Scolaro2022DevelopmentSolvers}
A.~Scolaro, C.~Fiorina, I.~Clifford, A.~Pautz, {Development of a semi-implicit
  contact methodology for finite volume stress solvers}, International Journal
  for Numerical Methods in Engineering 123~(2) (2022) 309--338.
\newblock \href {https://doi.org/https://doi.org/10.1002/nme.6857}
  {\path{doi:https://doi.org/10.1002/nme.6857}}.

\bibitem{Xing2018AProblems}
W.~Xing, C.~Song, F.~Tin-Loi, {A scaled boundary finite element based
  node-to-node scheme for 2D frictional contact problems}, Computer Methods in
  Applied Mechanics and Engineering 333 (2018) 114--146.
\newblock \href {https://doi.org/10.1016/J.CMA.2018.01.012}
  {\path{doi:10.1016/J.CMA.2018.01.012}}.

\bibitem{Aldakheel2020CurvilinearMechanics}
F.~Aldakheel, B.~Hudobivnik, E.~Artioli, L.~Beir{\~{a}}o~da Veiga, P.~Wriggers,
  {Curvilinear virtual elements for contact mechanics}, Computer Methods in
  Applied Mechanics and Engineering 372 (2020) 113394.
\newblock \href {https://doi.org/10.1016/J.CMA.2020.113394}
  {\path{doi:10.1016/J.CMA.2020.113394}}.

\bibitem{Popp2009AStrategy}
A.~Popp, M.~W. Gee, W.~A. Wall, {A finite deformation mortar contact
  formulation using a primal–dual active set strategy}, International Journal
  for Numerical Methods in Engineering 79~(11) (2009) 1354--1391.
\newblock \href {https://doi.org/https://doi.org/10.1002/nme.2614}
  {\path{doi:https://doi.org/10.1002/nme.2614}}.

\bibitem{DeLorenzis2012AMethod}
L.~De~Lorenzis, P.~Wriggers, G.~Zavarise, {A mortar formulation for 3D large
  deformation contact using NURBS-based isogeometric analysis and the augmented
  Lagrangian method}, Computational Mechanics 49 (2012) 1--20.
\newblock \href {https://doi.org/10.1007/s00466-011-0623-4}
  {\path{doi:10.1007/s00466-011-0623-4}}.

\bibitem{PavanaChand1998RemeshingProblems}
C.~PavanaChand, R.~KrishnaKumar, {Remeshing issues in the finite element
  analysis of metal forming problems}, Journal of Materials Processing
  Technology 75~(1) (1998) 63--74.
\newblock \href {https://doi.org/https://doi.org/10.1016/S0924-0136(97)00293-8}
  {\path{doi:https://doi.org/10.1016/S0924-0136(97)00293-8}}.

\bibitem{Sulsky1994AMaterials}
D.~Sulsky, Z.~Chen, H.~L. Schreyer, {A particle method for history-dependent
  materials}, Computer Methods in Applied Mechanics and Engineering 118~(1-2)
  (1994) 179--196.
\newblock \href {https://doi.org/10.1016/0045-7825(94)90112-0}
  {\path{doi:10.1016/0045-7825(94)90112-0}}.

\bibitem{Wang2019OnDistortion}
L.~Wang, W.~M. Coombs, C.~E. Augarde, M.~Cortis, T.~J. Charlton, M.~J. Brown,
  J.~Knappett, A.~Brennan, C.~Davidson, D.~Richards, A.~Blake, {On the use of
  domain-based material point methods for problems involving large distortion},
  Computer Methods in Applied Mechanics and Engineering 355 (2019) 1003--1025.
\newblock \href {https://doi.org/10.1016/J.CMA.2019.07.011}
  {\path{doi:10.1016/J.CMA.2019.07.011}}.

\bibitem{Coombs2020OnElasto-plasticity}
W.~M. Coombs, C.~E. Augarde, A.~J. Brennan, M.~J. Brown, T.~J. Charlton, J.~A.
  Knappett, Y.~Ghaffari~Motlagh, L.~Wang, {On Lagrangian mechanics and the
  implicit material point method for large deformation elasto-plasticity},
  Computer Methods in Applied Mechanics and Engineering 358 (2020) 112622.
\newblock \href {https://doi.org/10.1016/J.CMA.2019.112622}
  {\path{doi:10.1016/J.CMA.2019.112622}}.

\bibitem{Bardenhagen2001AnMaterial}
S.~G. Bardenhagen, J.~E. Guilkey, K.~M. Roessig, J.~U. Brackbill, W.~M. Witzel,
  J.~C. Foster, {An Improved Contact Algorithm for the Material Point Method
  and Application to Stress Propagation in Granular Material}, Computer
  Modeling in Engineering {\&} Sciences 2~(4) (2001) 509--522.
\newblock \href {https://doi.org/10.3970/cmes.2001.002.509}
  {\path{doi:10.3970/cmes.2001.002.509}}.

\bibitem{Homel2017FieldgradientMethod}
M.~A. Homel, E.~B. Herbold, {Field‐gradient partitioning for fracture and
  frictional contact in the material point method}, International Journal for
  Numerical Methods in Engineering 109~(7) (2017) 1013--1044.
\newblock \href {https://doi.org/10.1002/nme.5317}
  {\path{doi:10.1002/nme.5317}}.

\bibitem{Kakouris2019Phase-FieldEnergy}
E.~G. Kakouris, S.~P. Triantafyllou, {Phase-Field Material Point Method for
  dynamic brittle fracture with isotropic and anisotropic surface energy},
  Computer Methods in Applied Mechanics and Engineering 357 (2019).
\newblock \href {https://doi.org/10.1016/j.cma.2019.06.014}
  {\path{doi:10.1016/j.cma.2019.06.014}}.

\bibitem{Moutsanidis2019ModelingField}
G.~Moutsanidis, D.~Kamensky, D.~Z. Zhang, Y.~Bazilevs, C.~C. Long, {Modeling
  strong discontinuities in the material point method using a single velocity
  field}, Computer Methods in Applied Mechanics and Engineering 345 (2019)
  584--601.
\newblock \href {https://doi.org/10.1016/j.cma.2018.11.005}
  {\path{doi:10.1016/j.cma.2018.11.005}}.

\bibitem{Nairn2020NewFiltering}
J.~A. Nairn, C.~C. Hammerquist, G.~D. Smith, {New material point method contact
  algorithms for improved accuracy, large-deformation problems, and proper
  null-space filtering}, Computer Methods in Applied Mechanics and Engineering
  362 (4 2020).
\newblock \href {https://doi.org/10.1016/J.CMA.2020.112859}
  {\path{doi:10.1016/J.CMA.2020.112859}}.

\bibitem{Xiao2021DP-MPM:Fragmentation}
M.~Xiao, C.~Liu, W.~C. Sun, {DP-MPM: Domain partitioning material point method
  for evolving multi-body thermal–mechanical contacts during dynamic fracture
  and fragmentation}, Computer Methods in Applied Mechanics and Engineering 385
  (2021) 114063.
\newblock \href {https://doi.org/10.1016/J.CMA.2021.114063}
  {\path{doi:10.1016/J.CMA.2021.114063}}.

\bibitem{deVaucorbeil2021ModellingMethod}
A.~de~Vaucorbeil, V.~P. Nguyen, {Modelling contacts with a total Lagrangian
  material point method}, Computer Methods in Applied Mechanics and Engineering
  373 (1 2021).
\newblock \href {https://doi.org/10.1016/J.CMA.2020.113503}
  {\path{doi:10.1016/J.CMA.2020.113503}}.

\bibitem{Guilkey2021AMethod}
J.~Guilkey, R.~Lander, L.~Bonnell, {A hybrid penalty and grid based contact
  method for the Material Point Method}, Computer Methods in Applied Mechanics
  and Engineering 379 (2021) 113739.
\newblock \href {https://doi.org/10.1016/J.CMA.2021.113739}
  {\path{doi:10.1016/J.CMA.2021.113739}}.

\bibitem{Yamaguchi2021ExtendedMethod}
Y.~Yamaguchi, S.~Moriguchi, K.~Terada, {Extended B-spline-based implicit
  material point method}, International Journal for Numerical Methods in
  Engineering 122~(7) (2021) 1746--1769.
\newblock \href {https://doi.org/https://doi.org/10.1002/nme.6598}
  {\path{doi:https://doi.org/10.1002/nme.6598}}.

\bibitem{Hamad2017AMPM}
F.~Hamad, S.~Giridharan, C.~Moormann, {A Penalty Function Method for Modelling
  Frictional Contact in MPM}, Procedia Engineering 175 (2017) 116--123.
\newblock \href {https://doi.org/10.1016/J.PROENG.2017.01.038}
  {\path{doi:10.1016/J.PROENG.2017.01.038}}.

\bibitem{Huang2011ContactSimulation}
P.~Huang, X.~Zhang, S.~Ma, X.~Huang, {Contact algorithms for the material point
  method in impact and penetration simulation}, International Journal for
  Numerical Methods in Engineering 85~(4) (2011) 498--517.
\newblock \href {https://doi.org/10.1002/nme.2981}
  {\path{doi:10.1002/nme.2981}}.

\bibitem{Liu2020ILS-MPM:Particles}
C.~Liu, W.~C. Sun, {ILS-MPM: An implicit level-set-based material point method
  for frictional particulate contact mechanics of deformable particles},
  Computer Methods in Applied Mechanics and Engineering 369 (2020) 113168.
\newblock \href {https://doi.org/10.1016/J.CMA.2020.113168}
  {\path{doi:10.1016/J.CMA.2020.113168}}.

\bibitem{Guilkey2023CohesiveMethod}
J.~Guilkey, O.~Alsolaiman, R.~Lander, L.~Bonnell, J.~Cook, {Cohesive zones to
  model bonding in granular material with the material point method}, Computer
  Methods in Applied Mechanics and Engineering 415 (2023) 116260.
\newblock \href {https://doi.org/10.1016/J.CMA.2023.116260}
  {\path{doi:10.1016/J.CMA.2023.116260}}.

\bibitem{Chen2023DEM-enrichedMethod}
H.~Chen, S.~Zhao, J.~Zhao, X.~Zhou, {DEM-enriched contact approach for material
  point method}, Computer Methods in Applied Mechanics and Engineering 404
  (2023) 115814.
\newblock \href {https://doi.org/10.1016/J.CMA.2022.115814}
  {\path{doi:10.1016/J.CMA.2022.115814}}.

\bibitem{Bardenhagen2004TheMethod}
S.~Bardenhagen, E.~Kober, {The Generalized Interpolation Material Point
  Method}, Cmes-computer Modeling in Engineering {\&} Sciences (2004).
\newblock \href {https://doi.org/10.3970/CMES.2004.005.477}
  {\path{doi:10.3970/CMES.2004.005.477}}.

\bibitem{Charlton2017IGIMP:Deformations}
T.~J. Charlton, W.~M. Coombs, C.~E. Augarde, {iGIMP: An implicit generalised
  interpolation material point method for large deformations}, Computers and
  Structures 190 (2017) 108--125.
\newblock \href {https://doi.org/10.1016/J.COMPSTRUC.2017.05.004}
  {\path{doi:10.1016/J.COMPSTRUC.2017.05.004}}.

\bibitem{Sadeghirad2011ADeformations}
A.~Sadeghirad, R.~M. Brannon, J.~Burghardt, {A convected particle domain
  interpolation technique to extend applicability of the material point method
  for problems involving massive deformations}, International Journal for
  Numerical Methods in Engineering 86~(12) (2011) 1435--1456.
\newblock \href {https://doi.org/10.1002/NME.3110}
  {\path{doi:10.1002/NME.3110}}.

\bibitem{Tran2020AMethod}
Q.~A. Tran, W.~So{\l}owski, M.~Berzins, J.~Guilkey, {A convected particle least
  square interpolation material point method}, International Journal for
  Numerical Methods in Engineering 121~(6) (2020) 1068--1100.
\newblock \href {https://doi.org/10.1002/NME.6257}
  {\path{doi:10.1002/NME.6257}}.

\bibitem{Sadeghirad2013Second-orderInterfaces}
A.~Sadeghirad, R.~M. Brannon, J.~E. Guilkey, {Second-order convected particle
  domain interpolation (CPDI2) with enrichment for weak discontinuities at
  material interfaces}, International Journal for Numerical Methods in
  Engineering 95~(11) (2013) 928--952.
\newblock \href {https://doi.org/10.1002/NME.4526}
  {\path{doi:10.1002/NME.4526}}.

\bibitem{Wilson2021DistillationRemedy}
P.~Wilson, R.~W{\"{u}}chner, D.~Fernando, {Distillation of the material point
  method cell crossing error leading to a novel quadrature-based C0 remedy},
  International Journal for Numerical Methods in Engineering 122~(6) (2021)
  1513--1537.
\newblock \href {https://doi.org/10.1002/NME.6588}
  {\path{doi:10.1002/NME.6588}}.

\bibitem{Liang2019AnMethod}
Y.~Liang, X.~Zhang, Y.~Liu, {An efficient staggered grid material point
  method}, Computer Methods in Applied Mechanics and Engineering 352 (2019)
  85--109.
\newblock \href {https://doi.org/10.1016/J.CMA.2019.04.024}
  {\path{doi:10.1016/J.CMA.2019.04.024}}.

\bibitem{Nguyen2023TheApplications}
V.~P. Nguyen, A.~d. Vaucorbeil, S.~Bordas, {The Material Point Method: Theory,
  Implementations and Applications}, Scientific Computation, Springer Cham,
  2023.
\newblock \href {https://doi.org/10.1007/978-3-031-24070-6}
  {\path{doi:10.1007/978-3-031-24070-6}}.

\bibitem{Wang2024AnProblems}
M.~Wang, S.~Li, H.~Zhou, X.~Wang, K.~Peng, C.~Yuan, J.~Li, {An improved
  convected particle domain interpolation material point method for large
  deformation geotechnical problems}, International Journal for Numerical
  Methods in Engineering 125~(4) (2024) e7389.
\newblock \href {https://doi.org/10.1002/NME.7389}
  {\path{doi:10.1002/NME.7389}}.

\bibitem{Bing2019B-splineMethod}
Y.~Bing, M.~Cortis, T.~J. Charlton, W.~M. Coombs, C.~E. Augarde, {B-spline
  based boundary conditions in the material point method}, Computers {\&}
  Structures 212 (2019) 257--274.
\newblock \href {https://doi.org/10.1016/J.COMPSTRUC.2018.11.003}
  {\path{doi:10.1016/J.COMPSTRUC.2018.11.003}}.

\bibitem{Moutsanidis2020IGA-MPM:Method}
G.~Moutsanidis, C.~C. Long, Y.~Bazilevs, {IGA-MPM: The Isogeometric Material
  Point Method}, Computer Methods in Applied Mechanics and Engineering 372
  (2020) 113346.
\newblock \href {https://doi.org/10.1016/j.cma.2020.113346}
  {\path{doi:10.1016/j.cma.2020.113346}}.

\bibitem{Sulsky1995ApplicationMechanics}
D.~Sulsky, S.~J. Zhou, H.~L. Schreyer, {Application of a particle-in-cell
  method to solid mechanics}, Computer Physics Communications 87~(1-2) (1995)
  236--252.
\newblock \href {https://doi.org/10.1016/0010-4655(94)00170-7}
  {\path{doi:10.1016/0010-4655(94)00170-7}}.

\bibitem{Hollig2002WeightedProblems}
K.~H{\"{o}}llig, U.~Reif, J.~Wipper, {Weighted extended B-spline approximation
  of Dirichlet problems}, SIAM Journal on Numerical Analysis 39~(2) (2002)
  442--462.
\newblock \href {https://doi.org/10.1137/S0036142900373208}
  {\path{doi:10.1137/S0036142900373208}}.

\bibitem{Wriggers2006ComputationalMechanics}
P.~Wriggers, {Computational contact mechanics}, Computational Contact Mechanics
  (2006) 1--518\href {https://doi.org/10.1007/978-3-540-32609-0}
  {\path{doi:10.1007/978-3-540-32609-0}}.

\bibitem{Bonet2008NonlinearAnalysis}
J.~Bonet, R.~D. Wood, {Nonlinear Continuum Mechanics for Finite Element
  Analysis}, 2nd Edition, Cambridge University Press, Cambridge, 2008.
\newblock \href {https://doi.org/DOI: 10.1017/CBO9780511755446} {\path{doi:DOI:
  10.1017/CBO9780511755446}}.

\bibitem{DeLorenzis2011AAnalysis}
L.~De~Lorenzis, I.~Temizer, P.~Wriggers, G.~Zavarise, {A large deformation
  frictional contact formulation using NURBS-based isogeometric analysis},
  International Journal for Numerical Methods in Engineering 87~(13) (2011)
  1278--1300.
\newblock \href {https://doi.org/10.1002/NME.3159}
  {\path{doi:10.1002/NME.3159}}.

\bibitem{Hughes2005IsogeometricRefinement}
T.~Hughes, J.~Cottrell, Y.~Bazilevs, {Isogeometric analysis: CAD, finite
  elements, NURBS, exact geometry and mesh refinement}, Computer Methods in
  Applied Mechanics and Engineering 194~(39-41) (2005) 4135--4195.
\newblock \href {https://doi.org/10.1016/J.CMA.2004.10.008}
  {\path{doi:10.1016/J.CMA.2004.10.008}}.

\bibitem{Kakouris2018MaterialApproach}
E.~G. Kakouris, S.~P. Triantafyllou, {Material point method for crack
  propagation in anisotropic media: a phase field approach}, Archive of Applied
  Mechanics 88 (2018) 287--316.
\newblock \href {https://doi.org/10.1007/s00419-017-1272-7}
  {\path{doi:10.1007/s00419-017-1272-7}}.

\bibitem{Hammerquist2017AStability}
C.~C. Hammerquist, J.~A. Nairn, {A new method for material point method
  particle updates that reduces noise and enhances stability}, Computer Methods
  in Applied Mechanics and Engineering 318 (2017) 724--738.
\newblock \href {https://doi.org/10.1016/J.CMA.2017.01.035}
  {\path{doi:10.1016/J.CMA.2017.01.035}}.

\bibitem{Jiang2017AnMethod}
C.~Jiang, C.~Schroeder, J.~Teran, {An angular momentum conserving
  affine-particle-in-cell method}, Journal of Computational Physics 338 (2017)
  137--164.
\newblock \href {https://doi.org/10.1016/J.JCP.2017.02.050}
  {\path{doi:10.1016/J.JCP.2017.02.050}}.

\bibitem{Pretti2023AAnalysis}
G.~Pretti, W.~M. Coombs, C.~E. Augarde, B.~Sims, M.~Marchena~Puigvert, J.~A.~R.
  Guti{\'{e}}rrez, {A conservation law consistent updated Lagrangian material
  point method for dynamic analysis}, Journal of Computational Physics 485
  (2023) 112075.
\newblock \href {https://doi.org/10.1016/J.JCP.2023.112075}
  {\path{doi:10.1016/J.JCP.2023.112075}}.

\bibitem{Kikuchi1988ContactMethods}
N.~Kikuchi, J.~T. Oden, {Contact problems in elasticity: a study of variational
  inequalities and finite element methods}, SIAM, 1988.
\newblock \href {https://doi.org/10.1137/1.9781611970845}
  {\path{doi:10.1137/1.9781611970845}}.

\bibitem{Leichner2019ARepresentation}
A.~Leichner, H.~Andr{\"{a}}, B.~Simeon, {A contact algorithm for voxel-based
  meshes using an implicit boundary representation}, Computer Methods in
  Applied Mechanics and Engineering 352 (2019) 276--299.
\newblock \href {https://doi.org/10.1016/J.CMA.2019.04.008}
  {\path{doi:10.1016/J.CMA.2019.04.008}}.

\bibitem{Zhang2011MaterialFunction}
D.~Z. Zhang, X.~Ma, P.~T. Giguere, {Material point method enhanced by modified
  gradient of shape function}, Journal of Computational Physics 230~(16) (2011)
  6379--6398.
\newblock \href {https://doi.org/10.1016/J.JCP.2011.04.032}
  {\path{doi:10.1016/J.JCP.2011.04.032}}.

\bibitem{Johnson1985ContactMechanics}
K.~L. Johnson, {Contact Mechanics}, Cambridge University Press, Cambridge,
  1985.
\newblock \href {https://doi.org/DOI: 10.1017/CBO9781139171731} {\path{doi:DOI:
  10.1017/CBO9781139171731}}.

\end{thebibliography}


\begin{appendices}

\section{Auxiliary vectors for contact enforcement}\label{App:Contact}

In the two-dimensional case, matrices $\mathbf{C}^{nor}$ and $\mathbf{C}^{tan}$ are defined as 
\begin{equation}
\prescript{(t)}{}{\mathbf{C}}^{nor}_{} = -\left( \prescript{(t)}{}{ \mathbf{N} }^{  }_{  } - \frac{\prescript{(t)}{  }{ \textsl{g} }^{nor}_{}}{\prescript{(t)}{  }{ l }^{}_{ms}} \prescript{(t)}{  }{ \mathbf{P} }^{  }_{  } \right)
\label{eqn:C_nor}
\end{equation}
and
\begin{equation}
\prescript{(t)}{}{\mathbf{C}}^{tan}_{} = \prescript{(t)}{}{ \mathbf{T} }^{  }_{  } - \frac{\prescript{(t)}{  }{ \textsl{g} }^{tan}_{}}{\prescript{(t)}{  }{ l }^{}_{ms}} \prescript{(t)}{  }{ \mathbf{Q} }^{  }_{  },
\label{eqn:C_tan}
\end{equation}
respectively. Additionally, the $\mathbf{N}$, $\mathbf{Q}$, $\mathbf{T}$ and $\mathbf{P}$ matrices assume the following form
\begin{subnumcases}{}
\prescript{(t)}{}{\mathbf{N}}^{}_{} = \begin{bmatrix}
\prescript{(t)}{}{\mathbf{e}}^{nor}_{} \\
-(1-\prescript{(t)}{  }{ \beta }^{  }_{  } )\prescript{(t)}{}{\mathbf{e}}^{nor}_{} \\
-\prescript{(t)}{  }{ \beta }^{  }_{  }\prescript{(t)}{}{\mathbf{e}}^{nor}_{}
\end{bmatrix}
\label{eqn:N_2d} \\
\prescript{(t)}{  }{ \mathbf{Q} }^{  }_{  } = \begin{bmatrix}
                \mathbf{0} \\ 
                -\prescript{(t)}{  }{\mathbf{e} }^{ tan }_{  }  \\
                 \prescript{(t)}{  }{ \mathbf{e} }^{ tan }_{  } 
              \end{bmatrix} \label{eqn:Q_2d} \\
\prescript{(t)}{}{\mathbf{T}}^{}_{} = \begin{bmatrix}
\prescript{(t)}{}{\mathbf{e}}^{tan}_{} \\
-(1-\prescript{(t)}{  }{ \beta }^{  }_{  } )\prescript{(t)}{}{\mathbf{e}}^{tan}_{} \\
-\prescript{(t)}{  }{ \beta }^{  }_{  }\prescript{(t)}{}{\mathbf{e}}^{tan}_{}
\end{bmatrix}
\label{eqn:T_2d} \\
\prescript{(t)}{  }{ \mathbf{P} }^{  }_{  } = \begin{bmatrix}
                \mathbf{0} \\ 
                -\prescript{(t)}{  }{\mathbf{e} }^{ nor }_{  }  \\
                 \prescript{(t)}{  }{ \mathbf{e} }^{ nor }_{  } 
              \end{bmatrix} \label{eqn:P_2d}
\end{subnumcases}
with
\begin{equation}
\prescript{(t)}{}{\mathbf{u}}^{}_{} = \begin{bmatrix} 
\prescript{(t)}{}{\mathbf{u}}^{}_{s} \\
\prescript{(t)}{}{\mathbf{u}}^{}_{1} \\
\prescript{(t)}{}{\mathbf{u}}^{}_{2}
\end{bmatrix}.
\label{eqn:u_2d}
\end{equation}

\EQS{eqn:C_nor}-\eqref{eqn:u_2d} correspond to a discretised domain boundary formed by linear segments as explained in Section \ref{sec:boundary_form}. The matrices $\mathbf{C}^{nor}$ and $\mathbf{C}^{tan}$ for the one-dimensional case can be derived in a similar manner.

\end{appendices}

\end{document}